\documentstyle[12pt]{article}

\newtheorem{theorem}{Theorem}

\newtheorem{prop}{Proposition}
\newtheorem{lemma}{Lemma}
\newtheorem{cor}{Corollary}

\newcommand{\eqa}{\begin{eqnarray}}
\newcommand{\eeqa}{\end{eqnarray}}
\newcommand{\beq}{\begin{equation}}
\newcommand{\eeq}{\end{equation}}
\newcommand{\nn}{\nonumber}
\newcommand{\pal}{\partial}

\newcommand{\al}{\alpha}

\newcommand{\pf}{\noindent{\it Proof \ }}
\newcommand{\epl}{Lemma is proved.$\quad \Box$}
\newcommand{\ept}{Theorem is proved.$\quad \Box$}
\newcommand{\epp}{Proposition is proved.$\quad \Box$}
\newcommand{\wo}{\tilde\Omega}
\newcommand{\uu}{{\cal U}}
\newcommand{\cc}{{\cal C}}

\begin{document}


\begin{center}
{\LARGE Frobenius Manifolds And Virasoro Constraints}
\vskip 0.6cm
{Boris Dubrovin${}^*$ \ \ Youjin Zhang${}^{**}$}\\
\vskip 0.4cm
{\small {\it ${}^{*}$  SISSA, Via Beirut 2--4, 34014 Trieste, Italy}}\\
{\small {\it email: dubrovin@sissa.it}}\\
{\small {\it ${}^{**}$  Division of Mathematics, Graduate School of Science}}\\
{\small {\it Kyoto University, Kyoto 606-5802, Japan}}\\
{\small {\it email: youjin@kusm.kyoto-u.ac.jp}}
\end{center}
\begin{center}
\begin{minipage}{120mm}
\vskip 0.5 in
\begin{center}{\bf Abstract} \end{center}
{For an arbitrary Frobenius manifold a system of Virasoro
constraints is constructed. In the semisimple case these constraints
are proved to hold true in the genus one approximation. Particularly,
the genus $\leq 1$ Virasoro conjecture of T.Eguchi, K.Hori, M.Jinzenji,
and C.-S.Xiong and  of S.Katz is proved for smooth projective varieties
having semisimple
quantum cohomology.}
\end{minipage}
\end{center}
\vskip 0.5 in   

\section{Introduction}\par 

In a 2D topological field theory (TFT) one is
interested first of all in computation of correlators of the chiral
primary fields $\phi_1=1$, $\phi_2$, \dots, $\phi_n$ and of their
gravitational descendents $\tau_0(\phi_\alpha)=\phi_\alpha$,
$\tau_1(\phi_\alpha)$, $\tau_2(\phi_\alpha)$, \dots, $\alpha = 1, \dots ,
n$. Due to topological invariance these correlators are just numbers
$<\tau_p(\phi_\alpha)\, \tau_q(\phi_\beta)\dots >_g$ depending only on the
labels $(\alpha,p)$, $(\beta,q)$, \dots of the fields and on the genus 
$g$ of a 2D surface where the fields live. 
The free energy of the
theory is the generating function of these numbers
\beq\label{su4}
{\cal F}(T) = \sum_{g=0}^\infty {\cal F}_g (T)
\eeq
where $T$ is the infinite vector of coupling constants $T^{\alpha,p}$,
$\alpha=1, \dots , n$, $p=0, \, 1, \, 2, \dots$ and the genus $g$ part
is defined by
\beq
{\cal F}_g (T) =\big< e^{\sum_{\alpha =1}^n \sum_{p=0}^\infty
\tau_p(\phi_\alpha) T^{\alpha,p}}\big >_g .
\eeq

The correlators of the fields for different genera $g$ satisfy complicated
identities by now completely settled only for the simplest case of pure
topological gravity (only one primary field $\phi_1=1$). Of course,
the shape of these identities and their derivations depend on the concrete
model of 2D TFT. There are, however, some universal identities
that can probably be used as the defining relations
of some class of  2D TFTs.

Physicists construct 2D TFTs by twisting N=2 supersymmetric theories
\cite{Li, Witten1}.
The main consequence of this general approach \cite{DVV}
is a remarkable system
of differential equations for the function
\beq\label{su1}
F(t^1, \dots, t^n) := {\cal F}_0 (T) \big |_{T^{\alpha,p}=0 ~{\rm for}~
p>0, ~T^{\alpha,0}=t^\alpha, \, \alpha=1, \dots, n}.
\eeq
This  is WDVV associativity equations. They were also derived by
mathematicians \cite{KM, MS, RT}
 who were mainly concerned with the topological
sigma-models
known in the mathematical literature under the name {\it quantum
cohomology}. The mathematical derivations were based on a careful study
of intersection theory on the moduli spaces of Riemann surfaces with
markings and on the moduli spaces of stable maps of these surfaces
to smooth projective varieties. (In this case it is more convenient
to consider the correlators as elements of the so-called Novikov ring,
see Example 2.2 below.)  
This method proved to be successful
also in the analysis of the structure of genus one free energy \cite{DW,
Ge, Gi2}. 

The coordinate-free form of WDVV is the notion of Frobenius manifold
(see Sect. 2 below) proposed in \cite{D5}. 
The main motivation for the study of
Frobenius manifolds in \cite{D5}
 was the idea to construct all the building of
a given 2D TFT starting from the corresponding Frobenius manifolds,
i.e., starting from a solution (\ref{su1}) of WDVV.
Other motivations arrived when Frobenius manifolds appeared in the theory
of Gromov - Witten invariants, in K.Saito's theory of primitive forms
for isolated singularities, in geometry of invariants of reflection groups
and of their extensions etc. (see details in \cite{D3,D4,D6,DZ1,Hit,Man}). 
Physically, all
these  differently looking mathematical objects are just different models
of 2D TFT. Our idea is that, on the mathematical side, the theory of 
Frobenius manifolds is the common denominator unifying these rather
distant mathematical theories.

Remarkably, the original programme of reconstructing a 2D TFT (i.e., in
our setting, of the full free energy ${\cal F}(T)$) starting from a
solution (\ref{su1}) 
to WDVV seems to become realistic under the assumption
of semisimplicity of the Frobenius manifold (physically semisimple are
the models of 2D TFT with all massive perturbations \cite{CV}). 
Namely, there
is a universal procedure to reconstruct the components ${\cal F}_0(T)$and
${\cal F}_1(T)$ starting from an arbitrary semisimple Frobenius manifold.
Particularly, this procedure reproduces all topological recursion
relations of \cite{DW, Ge, Ho} for correlators in topological sigma-models,
i.e., the relations between the intersection numbers of Gromov - Witten
and Mumford - Morita - Miller cycles on the moduli spaces of stable maps
of punctured Riemann surfaces of the genus $g\leq 1$ to smooth projective
varieties. For genus 0 the reconstruction procedure of ${\cal F}_0(T)$ was
found in \cite{D1,D3}.
 Reconstruction of the genus 1 component ${\cal F}_1(T)$ for
an arbitrary semisimple Frobenius manifold was obtained in our paper 
\cite{DZ2}.

In the present paper we suggest new evidence supporting our vision of
semisimple Frobenius manifolds as the mathematical basement of the
building of 2D TFT. For an arbitrary Frobenius manifold we construct
a system of linear differential operators of the form
\eqa
&&L_m =\sum_{(\alpha,p), \, (\beta,q)} a_m^{\alpha,p; \, \beta,q}
{\partial^2
\over \partial T^{\alpha,p} \partial T^{\beta,q}}
+ \sum_{(\alpha,p), \, (\beta,q)} {b_m}_{\alpha,p}^ {\beta,q}\,
\tilde T^{\alpha,p}{\partial\over \partial T^{\beta,q}}
\nn\\
&&\quad
+\sum_{(\alpha,p), \, (\beta,q)} c^m_{\alpha,p; \, \beta,q}\,
\tilde T^{\alpha,p}\, \tilde T^{\beta,q} + {\rm const}\, \delta_{m,0},
\quad m\geq -1\label{su9}
\eeqa
with
\beq\label{su7}
\tilde T^{\alpha,p}=T^{\alpha,p}~~{\rm for} ~(\alpha,p)\neq (1,1), 
\quad
\tilde T^{1,1} = T^{1,1} -1
\eeq
such that, in the semisimple case, the partition function
\beq\label{su8}
Z(T) =e^{\sum_{g=0}^\infty{\cal F}_g(T)}
\eeq
satisfies, within $g\leq 1$ approximation, an infinite system of linear
equations
\beq\label{su2}
L_m Z(T) =0, ~~m\geq -1.
\eeq
Here $ a_m^{\alpha,p; \, \beta,q}$,
${b_m}_{\alpha,p}^ {\beta,q}$,
$c^m_{\alpha,p; \, \beta,q}$ are certain constants computed via the
monodromy data of the Frobenius manifold. Under certain nondegeneracy
assumption about the monodromy data the construction of the operators
$L_m$ can be extended to an arbitrary integer $m$. The operators $L_m$
satisfy the commutation relations of the Virasoro algebra
\beq\label{su3}
\left[ L_i, L_j\right] = (i-j) L_{i+j} +n {i^3-i\over 12} \delta_{i+j,0}.
\eeq
Because of this the equations (\ref{su2}) 
are called {\it Virasoro constraints}.

Our construction contains, as a particular case, the Virasoro constraints
in 2D topological gravity (in this case $n=1$, the Frobenius manifold is
trivial) established by M.Kontsevich in his proof \cite{Kon} of Witten's
conjecture \cite{Witten2}
 describing topological gravity by KdV hierarchy (see also
\cite{Witten3}). It
generalizes also the Virasoro constraints conjectured in \cite{Egu3} in the
framework of quantum cohomology. Recently the quantum cohomology Virasoro
constraints of \cite{Egu3} were proved \cite{LT} in the genus zero
approximation. However, the genus zero result does not suffice to fix
uniquely all the coefficients of the Virasoro operators and, therefore,
it does not suffice to prove the commutation relations (\ref{su3}).

One may ask where do the Virasoro constraints come from? In our opinion,
validity of the constraints proved for $g\leq 1$ in the general setting in
the present paper is a manifestation of a general involvment of integrable
bihamiltonian hierarchies of KdV-type into an arbitrary 2D TFT. Due to an
idea of \cite{Witten2},
 the unknown functions of the hierarchy are the following
correlators
\beq\label{su14}
t_\alpha(T) := \sum_{g=0}^\infty \left< \tau_0(\phi_\alpha) \tau_0(\phi_1)
e^{\sum_{\beta, \, q}\tau_q(\phi_\beta) T^{\beta,q}}\right>_g, ~~\alpha=1,
\dots, n
\eeq
as the functions of the ``times'' $T^{\beta,q}$. The differential
equations
for the correlators can be considered as a hierarchy of commuting
evolutionary PDEs
\beq\label{su6}
\partial_{T^{\alpha,p}}t=\partial_X K_{\alpha,p}(t, t_X, t_{XX}, \dots)
\eeq
on the loop space 
$$
L(M) = \{(t^1(X), \dots, t^n(X))\, |\, X\in S^1\}.
$$
Here $X=T^{1,0}$ is one of the couplings. The partition function (\ref{su8})
is the tau-function of a particular solution to the hierarchy specified
by the {\it string equation}
\beq
L_{-1} Z \equiv \sum T^{\alpha, p+1} {\partial Z\over \partial
T^{\alpha,p}} +{1\over 2} \eta_{\alpha\beta} T^{\alpha,0} T^{\beta,0} Z
-{\partial Z\over \partial X} =0.\label{su5}
\eeq
The Virasoro operators correspond to the symmetries of the hierarchy. The
operator $L_{-1}$ correspond to a Galilean-type symmetry; $L_0$
corresponds to the scaling symmetry; other symmetries can be derived using
the recursion procedure (see, e.g., \cite{Soliton} and references therein) 
due to bihamiltonian structure of the hierarchy.

This scheme works perfectly well in the particular case of topological
gravity where the hierarchy (\ref{su6}) 
coincides with KdV \cite{Kon}. For other models
of 2D TFT the hierarchy and the
bihamiltonian structure are to be constructed. 
(This was done for $g=0$
in \cite{D1,D3} and for $g=1$ in our paper \cite{DZ2}, see below.) 
So, the programme
of
reconstruction of a 2D TFT from a given Frobenius manifold contains
a subprogramme of realization of the following identification
$$
\mbox{semisimple Frobenius manifolds}\equiv\begin{array}{l}
\mbox{moduli of integrable bihamiltonian} \\
\mbox{hierarchies of KdV type.}
\end{array}
$$

The last important idea to be explained in this introduction is a
technical trick of introducing of a small parameter $\epsilon$ in the
genus expansions (\ref{su8}) and (\ref{su14}). 
In physical literature the small parameter
$\epsilon$ is called {\it the string coupling constant}. Introduction of the
parameter $\epsilon$ is based on the observation, essentially due to
\cite{DW},
that the genus $g$ part of the free energy is a homogeneous function
of the degree $2-2g$ of the shifted couplings (\ref{su7}). 
So, after the rescaling
\beq\label{su10}
\tilde T^{\alpha,p}\mapsto \epsilon \tilde T^{\alpha,p}
\eeq
the genus expansion (\ref{su8}) reads
\beq\label{su11}
Z(T;\epsilon) =\exp{\sum_{g=0}^\infty \epsilon^{2g-2} {\cal F}_g(T)}.
\eeq

This trick is very convenient when separating the contributions of
different genera. Observe that the partition function (\ref{su11}) 
is annihilated
by the operator
\beq
D:= \sum \tilde T^{\alpha,p} \partial _{T^{\alpha,p}}
+\epsilon\partial_\epsilon+\frac{n}{24}.
\eeq

From the point of view of integrable hierachies the genus expansion
corresponds to the small dispersion expansion of (\ref{su6}) 
(we also rescale
$X\mapsto \epsilon X$)
\beq\label{su15}
\partial_{T^{\alpha,p}} t =\partial_X K^{(0)}_{\alpha,p} (t) + \epsilon^2
\partial_X K^{(1)}_{\alpha,p} (t, t_X, t_{XX})+\dots.
\eeq
It turns out \cite{DZ2} 
that for an arbitrary semisimple Frobenius manifold 
the first two terms of the bihamiltonian hierarchy can be constructed
(see \cite{D1,D3} for $K^{(0)}_{\alpha,p}$ and \cite{DZ2} 
for $K^{(1)}_{\alpha,p}$). Moreover, they are determined uniquely by the
properties of topological correlators.

Applying the Virasoro operators (\ref{su9}), 
where one is to do the same rescaling
(\ref{su10}), to the partition function (\ref{su11}) one obtains a series
$$
L_m Z(T;\epsilon) = \left(\sum_{g=0}^\infty {\cal
A}_{m,g}\epsilon^{2g-2}\right) Z(T;\epsilon)
$$
where the coefficients ${\cal A}_{m,g}$ are expressed via ${\cal F}_0$, 
\dots, ${\cal F}_g$ and their derivatives. Particularly,
\eqa
&&{\cal A}_{m,0}=
\sum a_m^{\alpha,p;\, \beta,q}\, \partial_{T^{\alpha,p}}{\cal F}_0\,
\partial_{T^{\beta,q}}{\cal F}_0
+\sum {b_m}_{\alpha,p}^{\beta,q}\,\tilde T^{\alpha,p}\,
\partial_{T^{\beta,q}}{\cal F}_0 
\nn\\
&&\quad+\sum c^m_{\alpha,p;\, \beta,q}\, \tilde
T^{\alpha,p}\,\tilde T^{\beta,q},\label{su12}
\\
&&\nn\\
&&{\cal A}_{m,1}=\sum a_m^{\alpha,p;\, \beta,q}\left(
\partial_{T^{\alpha,p}}\partial_{T^{\beta,q}} {\cal F}_0
+ \partial_{T^{\alpha,p}} {\cal F}_1\, \partial_{T^{\beta,q}}{\cal F}_0
+ \partial_{T^{\alpha,p}} {\cal F}_0\, \partial_{T^{\beta,q}}{\cal
F}_1\right)
\nn\\
&&\quad+ \sum {b_m}_{\alpha,p}^{\beta,q}\, \tilde T^{\alpha,p}\,
\partial_{T^{\beta,q}} {\cal F}_1 + {\rm const} \delta_{m,0}.\label{su13}
\eeqa

The main result of this paper, besides the construction of the Virasoro
operators (\ref{su9}), is

\vskip 0.4cm
\noindent{\bf Main Theorem}\hskip 0.4cm {\it
1. For an arbitrary Frobenius manifold
the genus zero Virasoro constraints
\beq\label{su16}
{\cal A}_{m,0}=0, ~~m\geq -1
\eeq
hold true.

2. For an arbitrary semisimple Frobenius manifold the genus 1 Virasoro
constraints
\beq\label{su17}
{\cal A}_{m,1}=0, ~~m\geq -1
\eeq
hold true.}

\medskip
In the setting of quantum cohomology of a variety $X$ having $H^{\rm
odd}(X)=0$ the first part of our theorem gives a new proof of the main
result of \cite{LT}. The second part, particularly, proves the Virasoro
conjecture of Eguchi, Hori, Jinzenji and Xiong \cite{Egu3, Egu2} and of 
S. Katz in the genus one approximation for smooth projective varieties $X$
having semisimple quantum cohomology. According to the conjecture of
\cite{D6} the semisimplicity 
takes place for those Fano varieties $X$ on which
a full system of exceptional objects, in the sense of \cite{BP}, 
exists in the
derived category of
coherent sheaves $Der^b(Coh(X))$. In the cohomology of such  a variety $X$
only  $H^{k,k}(X)\neq 0$. Because of this we do not consider in the
present
paper the Virasoro constraints for Frobenius supermanifolds neither
the extension of the Virasoro algebra due to the presence of nonanalytic
cycles, i.e., those belonging to $H^{k,l}(X)$ with $k\neq l$. Here $H(X)=
\oplus H^{k,l}(X)$ is the Hodge decomposition.

Our theorem also
fixes uniquely all the
coefficients of the Virasoro operators conjectured in \cite{Egu3,Egu2}. 
As a byproduct
of the proof of the second part of the theorem we discover
a simple way of expressing elliptic Gromov - Witten invariants via
rational ones.

The structure of the paper is as follows. In Sect. 2 we present a
Sugawara-type (see \cite{Kac}) construction for the needed Virasoro operators.
In Sect. 3 we recall the necessary information about semisimple Frobenius
manifolds and the construction of ${\cal F}_0(T)$ and ${\cal F}_1(T)$.
In Sect. 4 we sketch the proof of the main theorem.
\smallskip

{\bf Acknowledgment.} We are grateful to C.-S.Xiong for a very instructive
discussion of the results of \cite{Egu2}. Y. Z. is supported by the
Japan Society for the Promotion of Science, he is grateful to 
M. Jimbo for advises and encouragements.  

\setcounter{equation}{0}
\vskip 0.5cm
\section{A family of representations of Virasoro algebra}

We describe first the parameters of the family. Let $V$ be $n$-dimensional
complex vector space equipped with a symmetric nondegenerate bilinear form
$<~,~>$. The parameters of our Virasoro algebras consist of:

\smallskip

1. An antisymmetric linear operator
\beq\label{jw15}
\mu: V\to V, ~~<\mu(a),b>+<a,\mu(b)>=0.
\eeq
For the sake of simplicity in this paper we will consider
only the case of diagonalizable operators $\mu$.
For any eigenvalue $\lambda \in {\rm spec}\, \mu$ denote $V_\lambda$
the subspace of $V$
 consisting of all eigenvectors with the eigenvalue $\lambda$.
Let $\pi_\lambda:V\to V_\lambda$ be the
projector,
\beq\label{ad1}
\pi_\lambda |_{V_\lambda} = {\rm id}, ~~ \pi_\lambda (V_{\lambda'})=0
~{\rm for}~ \lambda'\neq \lambda.
\eeq
For any linear operator $A:V\to V$ and for any integer $k$ denote
\beq\label{ad2}
A_k := \sum_{\lambda\in {\rm spec}\, \mu} \pi_{\lambda+k} A \pi_\lambda.
\eeq
By the construction
\beq\label{ad3}
A_k(V_\lambda) \subset V_{\lambda+k}.
\eeq

\smallskip
2. A {\it $\mu$-nilpotent} linear operator $R:V\to V$. By definition this
means that
\beq\label{ad4}
R = R_1 \oplus R_2 \oplus \dots
\eeq
and the components $R_k$ (defined as in (2.3)) for any $k$ satisfy 
\beq\label{ad41}
\left[ \mu, R_k\right] = k R_k
\eeq
and also the following symmetry conditions
\beq\label{ad5}
<R_k(a), b>+(-1)^k <a, R_k(b)> =0, ~~k=1, \, 2, \dots, .
\eeq
These can be recasted into the form
\beq\label{ad6}
\{ R(a), b\} + \{ a, R(b)\} =0 ~~{\rm for ~ any}~ a, \, b \in V
\eeq
where the bilinear form $\{ ~,~\}$ on $V$ is defined by
\beq\label{jw1}
\{ a,b\} := \left< e^{\pi i \mu} a,b\right>.
\eeq
\smallskip
\noindent
{\bf Remark 2.1.} The operator $R$ must be zero if the operator $\mu$ is
{\it nonresonant}. By definition this means that any two eigenvalues of
$\mu$ do not differ by a nonzero integer. In the opposite case the
operator $\mu$ is called {\it resonant}.

\vskip 0.4cm
We define an equivalence relation between two sets $\left( V^{(i)},
<~,~>^{(i)}, \mu^{(i)}, R^{(i)}\right)$, $i=1, \, 2$ of the parameters.
An equivalence is established by an isomorphism 
\beq\label{jw4-1}
G: V^{(1)}\to V^{(2)}
\eeq
of the linear spaces satisfying
$$
G=G_0 \oplus G_1 \oplus \dots, 
$$
$G_0$ is an isometry
$$
\left< G_0(a), G_0(b)\right>^{(2)} =\left< a, b\right> ^{(1)}
$$
satisfying
$$
G_0\mu^{(1)} = \mu^{(2)} G_0
$$
other components $G_k$ with $k>0$ satisfy
$$
\mu^{(2)} G_k -G_k \mu^{(1)} = k G_k
$$
$$
G\, R^{(1)} = R^{(2)}\,G
$$
$$
\left\{ G(a), G(b)\right\}^{(2)} = \{ a,b\}^{(1)} ~~{\rm for ~any }~a, \,
b\in V^{(1)}.
$$
The bilinear forms $\{ ~,~\}^{(1)}$ and $\{ ~,~\}^{(2)}$ on the spaces
$V^{(1)}$ and $V^{(2)}$ resp. are defined by (\ref{jw1}).
\smallskip

\noindent
{\bf Remark 2.2.} As was explained in \cite{D4}, the  equivalence class
of a quadruple $\left( V, <~,~>, \mu, R\right)$ is the set of monodromy
data
at $z=0$ of the linear differential equation
\beq\label{jw2}
z{dy\over dz} =A(z)y
\eeq
where the $n\times n$ matrix $A(z)$ satisfies
\beq
A(z) = \sum_{k\geq 0} A_k z^k, ~~A_0=\mu,\nn
\eeq
\beq\label{jw5}
A^T(-z)\eta + \eta A(z)=0
\eeq
for  constant symmetric nondegenerate matrix $\eta$. Hear $\eta$ is the
Gram matrix of the bilinear form $<~,~>$ w.r.t. a basis in $V$. We denote
the matrices of the operators $\mu$, $R$ w.r.t. this basis by the same
letters $\mu$, $R$.
The monodromy data,
according to the Birkhoff-style definition \cite{CodLev} are defined as the
equivalence class of the system (\ref{jw2}) w.r.t. gauge transformations
of the form
\beq\label{jw6}
y\mapsto \Theta(z) y, \quad \Theta(z)= \sum_{k\geq 0} \Theta_k z^k,
\quad
\Theta^T(-z) \eta \Theta(z) =\eta.
\eeq
Particularly, this means that the system (\ref{jw2}) 
posesses a fundamental matrix
of the form
\beq\label{jw7}
Y(z)=\Theta(z) z^\mu z^R
\eeq
where the matrices $\Theta(z)$, $\eta$, $\mu$, $R$ satisfy the above
conditions.
For two equivalent quadruples of the monodromy data the correspondent
matrices $\Theta^{(1)}(z)$ and $\Theta^{(2)}(z)$ are related by right
multiplication by the matrix $G_0 + z G_1 + z^2 G_2 +\dots$.

\vskip 0.4cm
We proceed now to the construction of the Virasoro algebra for any data
\newline
$\left( V, <~,~>, \mu, R\right)$. Choose a basis $e_1$, \dots, $e_n$ in
$V$. Denote
\beq
\eta_{\alpha\beta}=<e_\alpha,e_\beta>.
\eeq
Introduce the inverse matrix
\beq
\left( \eta^{\alpha\beta}\right) = \left( \eta_{\alpha\beta}\right)^{-1}.
\eeq
We start with introducing of a Heisenberg algebra of $\left( V,
<~,~>\right)$. It has generators ${\bf 1}$ and $a^\alpha_k$, $\alpha=1,
\dots,
n$, $k\in {\bf Z}$ obeying the commutation relations
\beq\label{jw16}
\left[{\bf 1}, a^\alpha_k\right] =0,
\quad
\left[ a^\alpha_k, a^\beta_l\right] =(-1)^k \eta^{\alpha\beta}
\delta_{k+l+1,0}\cdot {\bf 1}.
\eeq
The normal ordering is defined by
\beq\label{ad7}
: a^\alpha_k a^\beta_l: =\left\{
\begin{array}{ll}
a^\beta_l a^\alpha_k, &{\rm if}~l<0,
~k\geq 0\\
 a^\alpha_k a^\beta_l, &{\rm otherwise}.
\end{array}
\right.
\eeq
Introduce vector-valued operators
\beq
a_k =a^\alpha_k e_\alpha,~~k\in {\bf Z}.
\eeq
For any integer $m$ we define matrices in $End(V)$
\beq
P_m(\mu, R) :=\left\{
\begin{array}{cl}
\left[ e^{R\, \partial_x} \prod_{j=0}^m \left( x+\mu+j-{1\over 2}\right)
\right]_{x=0}, &m\geq 0\\
1, &m=-1\\
\left[ e^{R\, \partial_x} \prod_{j=1}^{-m-1} \left( x+\mu-j-{1\over
2}\right)^{-1}
\right]_{x=0}, &m< -1.
\end{array}\right.
\eeq

\vskip 0.3cm
\noindent
{\bf Remark 2.3.} The matrices $P_m(\mu,R)$ for $m<-1$ are defined only if
the spectrum of $\mu$ does not contain half-integers.

\vskip 0.3cm
Now we define the operators of the Virasoro algebra by a Sugawara-type
construction 
\eqa
L_m&=&{1\over 2} \sum_{k,\,l} (-1)^{k+1} : \left< a_l, \left[
P_m(\mu-k,R)\right]_{m-1-l-k} a_k\right> :
\nn\\
&& +{1\over 4} \delta_{m,0}\,  
{\rm tr}\left( {1\over 4} -\mu^2\right)\, {\bf 1}, \quad m\in {\bf Z}.
\label{jw3}
\eeqa
Here the components $[~~]_{q}$ for any integer $q$ are defined in 
(\ref{ad2}). 
Due to Remark 2.3 the operators $L_m$ with $m\geq -1$ are always defined,
another half, $L_{m<-1}$, is defined only under the assumption
${\rm spec}\, \mu \cap {1\over 2}+{\bf Z} =\emptyset$. 
\begin{theorem}
The operators (\ref{jw3}), when well-defined, satisfy
the Virasoro commutation relations (\ref{su3}) with the central charge $n$.
\end{theorem}

\pf For an arbitrary constant $z$ the following
identity is easily checked to hold true:
\beq\label{id1}
\eta\,[P_m(\mu+z+1,R)]_q=(-1)^{m+q+1}\,[P_m(\mu-m-z+q,R)]_q^T\,\eta.
\eeq
Using this identity we can express the $L_m$ operators in the form
\eqa
&&L_m=(1-\delta_{m,0})\,\sum_{0\le i\le m-1}\,\sum_{0\le k\le m-i-1}
\frac12\,(-1)^{k+1}\,\left<[P_m(\mu-k,R)]_i\,
a_k, a_{m-i-k-1}\right>
\nn\\
&&\quad+\sum_{k\ge 0}\,\sum_{0\le i\le m+k}(-1)^k
\left<[P_m(\mu+k+1,R)]_i\, a_{-k-1}, a_{m+k-i}\right>
\nn\\
&&\quad+\sum_{k,l\ge 0} \frac12\,(-1)^k \left<
[P_m(\mu+k+1,R)]_{l+k+m+1}\, a_{-k-1},
a_{-l-1}\right>
\nn\\
&&\quad+{1\over 4} \delta_{m,0}\,  
{\rm tr}\left( {1\over 4} -\mu^2\right)\, {\bf 1}, \quad m\ge 0,
\\
&& \nn\\
&&L_m=\sum_{i\ge 0}\,\sum_{k\ge i-m} (-1)^k\left<[P_m(\mu+k+1,R)]_i\,
 a_{-k-1}, a_{m+k-i}\right>
\nn\\
&&\quad+\sum_{i\ge 0}\,\sum_{0\le k\le i-m-1}
\frac12\,(-1)^k \left<[P_m(\mu+k+1,R)]_i\,
a_{-k-1}, a_{m+k-i}\right>,
\nn\\&&\quad m\le -1.
\eeqa
From the above expression of the $L_m$ operators and the
identity (\ref{id1}) it follows, for any $m, m'\ge 0, m+m'\ne 0$, that
the commutator $[L_m,L_{m'}]$ can be written as
\eqa
&&\{\sum_{0\le q\le m+m'-1}\ \sum_{0\le k\le m+m'-q-1}\ 
\sum_{i+i'=q}\ \frac12\,(-1)^{m+i+k+1}\times
\nn\\
&&\quad\times \left<[P_m(\mu-m-m'+q+k+1,R)]_i\,a_{m+m'-q-k-1},
[P_{m'}(\mu-k,R)]_{i'}\,a_k\right>
\nn\\
&&-\sum_{0\le q\le m+m'-1}\ \sum_{0\le k\le m+m'-q-1}\
\sum_{i+i'=q}\ \frac12\,(-1)^{m'+i'+k+1}\times
\nn\\
&&\quad\times \left<[P_{m'}(\mu-m-m'+q+k+1,R)]_{i'}\,a_{m+m'-q-k-1},
[P_{m}(\mu-k,R)]_{i}\,a_k\right>\}
\nn\\ 
&&\nn\\
&&+\{\sum_{k\ge 0}\ \sum_{0\le q\le m+m'+k}\ \sum_{i+i'=q}
\ (-1)^{m+i+k}\times
\nn\\
&&\quad\times \left<[P_{m'}(\mu+k+1,R)]_{i'}\,a_{-k-1},
[P_m(\mu-m-m'+q-k,R)]_i\,a_{m+m'-q+k}\right>
\nn\\
&&-\sum_{k\ge 0}\ \sum_{0\le q\le m+m'+k}\ \sum_{i+i'=q}
\ (-1)^{m'+i'+k}\times
\nn\\
&&\quad\times \left<[P_m(\mu+k+1,R)]_{i}\,a_{-k-1},
[P_{m'}(\mu-m-m'+q-k,R)]_{i'}\,a_{m+m'-q+k}\right>\}
\nn\\
&&\nn\\
&&+\{\sum_{k, k'\ge 0}\ \sum_{0\le i\le m+m'+k+k'+1}\ \frac12\,
(-1)^{m+k'+i}\times
\nn\\
&&\quad \times \left<[P_m(\mu+k+1,R)]_i\,a_{-k-1},
[P_{m'}(\mu+k'+1,R)]_{m+m'+k+k'-i+1}\,a_{-k'-1}\right>
\nn\\
&&-\sum_{k, k'\ge 0}\ \sum_{0\le i'\le m+m'+k+k'+1}\ \frac12\,
(-1)^{m'+k+i'}\times
\nn\\
&&\quad \times \left<[P_{m'}(\mu+k'+1,R)]_{i'}\,a_{-k'-1},
[P_{m}(\mu+k+1,R)]_{m+m'+k+k'-i'+1}\,a_{-k-1}\right>\}
\nn\\
&&\nn\\
&&-\frac12\,\sum_{0\le i\le m-1}\,\sum_{0\le k\le m-i-1}
{\rm tr}\left([P_m(\mu-k,R)]_i\,
[P_{m'}(\mu+m-i-k,R)]_{m+m'-i}\right)
\nn\\
&&+\frac12\,\sum_{0\le i'\le m'-1}\,\sum_{0\le k'\le m'-i'-1}
{\rm tr}\left([P_{m'}(\mu-k',R)]_{i'}\,
[P_{m}(\mu+m'-i'-k',R)]_{m+m'-i'}\right).\nn
\eeqa
The traces which appeared in the last two sums vanish 
(for $m+m'>0$) due to the fact that
\beq
[\mu,[R^k]_l]=l\,[R^k]_l,
\eeq
the remaining terms are then summed up to give
the desired $(m-m')\,L_{m+m'}$ due to the following identity
\eqa
&&\sum_{i+i'=q}(-1)^{m+i}\,[P_m(\mu-m-m'+q+z+1,R)]_i^T\,\eta\,
[P_{m'}(\mu-z,R)]_{i'}
\nn\\
&&-\sum_{i+i'=q}(-1)^{m'+i'}\,[P_{m'}(\mu-m-m'+q+z+1,R)]_{i'}^T\,\eta\,
[P_{m}(\mu-z,R)]_{i}
\nn\\
&&=(m-m')\,\eta\,[P_{m+m'}(\mu-z,R)]_q,\quad {\mbox {for}}\ m, m'\ge 0.
\eeqa

We can compute the commutator $[L_m,L_{m'}]$ for other cases of 
$m, m'\in {\bf Z}$ in a similar way. \ept 

\medskip
The algebra generated by the operators $L_m$ will be denoted by
$Vir\left( V,<~,~>, \mu, R\right)$.
It can be readily seen that two equivalent sets of data produce isomorphic
Virasoro algebras.

A particular realization of the Heisenberg algebra will be in the space
of functions of an infinite number of variables 
$T^{\alpha,p}$, $\alpha=1,
\dots, n$, $p=0,\, 1, \, 2, \dots$. We put
\beq\label{jw9}
a^\alpha_k =\left\{
\begin{array}{ll}
\eta^{\alpha\beta} {\partial\over \partial
T^{\beta,k}}, &k\geq 0\\
(-1)^{k+1} T^{\alpha, -k-1}, &k<0.
\end{array}\right.
\eeq
We also introduce vectors
$$
T^k =T^{\alpha,k} e_\alpha\in V\otimes {\bf C}[[T]], ~k\geq 0
$$
and covectors
$$
\partial_{T^k} = e^\alpha \partial_{T^{\alpha,k}}\in V^*\otimes Der
\left( {\bf C}[[T]]\right), ~k\geq 0,
$$
where $e^1$, \dots, $e^n\in V^*$ is the basis dual to $e_1$ \dots, $e_n\in 
V$. Denote by the same symbol $<~,~>$ the natural pairing $V\otimes V^*\to
{\bf C}$ and also the bilinear form on $V^*$ dual to that on $V$. Then the
first few Virasoro operators look as follows:
\eqa\label{jw10}
&&L_{-1}=\sum_{p\geq 1} \left< T^p ,\partial_{T^{p-1}}\right> +{1\over 2}
\left< T^0, T^0\right>
\\
&&L_0 = \sum _{p\geq 0} \left< \left( p+{1\over 2} +\mu\right) T^p,
\partial_{T^p}\right>
+\sum_{p\geq 1}\sum_{1\leq r\leq p} \left< R_r T^p,
\partial_{T^{p-r}}\right>
\nn\\
&&\quad
+{1\over 2} \sum_{p,\, q}(-1)^q\, \left< R_{p+q+1} T^p, T^q\right>
+{1\over 4} {\rm tr}\left( {1\over 4} -\mu^2\right).\label{jw11}
\eeqa
To write down the formulae for $L_1$ and $L_{\pm 2}$ we introduce the
matrices
$R_{k,l}\in End(V)$ putting
\eqa
R_{0,0} &=&1\nn\\
R_{k,0} &=&0, ~k>0\nn\\
R_{k,l} &=& \left[ R^l\right]_k=\sum_{i_1+\dots +i_l=k} R_{i_1} \dots
R_{i_l}, ~l>0.\nn
\eeqa
Then
\eqa
&&L_1 =\sum_{p\geq 0} \left< \left( p+{1\over 2} +\mu\right)\left(p+{3\over
2} + \mu\right) T^p, \partial_{T^{p+1}}\right>
\nn\\
&&\quad+\sum_{p\geq 0} \sum_{1\leq r\leq p+1} \left<
R_r \left( 2 p + 2 + 2 \mu\right) T^p, \partial_{T^{p-r+1}}\right>
\nn\\
&&\quad
+\sum_{p\geq 1} \sum_{2\leq r\leq p+1} \left< R_{r,2} T^p,
\partial_{T^{p-r+1}}\right>
+{1\over 2} \left< \partial_{T^0}\left({1\over 2}+\mu\right),
\partial_{T^0} \left( {1\over 2} +\mu\right)\right>
\nn\\
&&\quad
+\sum_{p,\, q\geq 0} (-1)^q \left< R_{p+q+2} \left( p+\mu+1\right) T^p,
T^q\right>
+{1\over 2} \sum_{p,\, q\geq 0} (-1)^q \left< R_{p+q+2,2} T^p, T^q\right>,
\label{jw12}
\\
&&\nn\\
&&
L_2=\sum_{p\geq 0} \left< \left( p+{1\over 2}+\mu\right) \left( p+{3\over
2} +\mu\right) \left(p+{5\over 2}+\mu\right) T^p,
\partial_{T^{p+2}}\right>
\nn\\
&&\quad
+\sum_{p\geq 0} \sum_{1\leq r\leq p+2}
\left< R_r \left[ 3\left( p+{1\over 2}+\mu\right)^2 + 6 \left( p+{1\over
2}+\mu\right) + 2\right] T^p , \partial_{T^{p-r+2}}\right>
\nn\\
&&\quad
+\sum_{p\geq 0}\sum_{2\leq r\leq p+2} \left< R_{r,2}\left(3 p +{9\over 2}
+ 3 \mu\right) T^p, \partial _{T^{p-r+2}}\right>
+\sum_{p\geq 1}\sum_{3\leq r\leq p+2} \left< R_{r,3} T^p,
\partial_{T^{p-r+2}}\right>
\nn\\
&&\quad
+\left< \partial_{T^1}\left({1\over 2}-\mu\right) ,
\partial_{T^0}\left({1\over 2}-\mu\right) \left({3\over 2}
-\mu\right)\right>
+{1\over 2}\left< \partial_{T^0} R_1, \partial_{T^0}\left({1\over 4} +
3\mu - 3\mu^2\right)\right>
\nn\\
&&\quad
+{1\over 2}\sum_{p, \, q\geq 0} (-1)^q
\left< \left( R_{p+q+3,2}+3 R_{p+q+3,2}\left( p+\mu+{3\over 2}\right)
\right.\right.
\nn\\
&&\quad
+\left.\left.
R_{p+q+3}\left[ {3\over 4} (2 p + 2 \mu + 3)^2 -1\right] \right) T^p, T^q
\right>,
\label{jw13}
\\
&&\quad \nn\\
&&
L_{-2}=\sum_{p\geq 2}\left<\left(p-{1\over 2} + \mu\right)^{-1} T^p,
\partial_{T^{p-2}}\right>
\nn\\
&&\quad
+\sum_{k\geq 1}\sum_{l\geq k} \sum_{p\geq l+2} (-1)^k
\left< R_{l,k}\left(p-{1\over 2}+\mu\right)^{-k-1} T^p,
\partial_{T^{p-l-2}}\right>
\nn\\
&&\quad
+\left< \left({1\over 2}-\mu\right)^{-1}  T^0,  T^1\right>
\nn\\
&&\quad
+\frac12\,\sum_{k\geq 1} \sum_{l\geq k} \sum_{0\leq p \leq l+1} (-1)^{l+p+k+1}
\left< R_{l,k} \left(p-{1\over 2}+\mu\right)^{-k-1} T^p, 
T^{l-p+1}\right>.\label{jw14}
\eeqa

\vskip 0.4cm
\noindent
{\bf Example 2.1.} For $n=1$ it must be $\mu=0$, $R=0$. So
$$
L_m={1\over 2} \sum_k (-1)^{k+1} P_m(-k) :a_k a_{m-1-k}:
+{1\over 16} \delta_{m,0}
$$
where 
$$
P_m(x)=\left\{
\begin{array}{cl}
\prod_{j=0}^m\left( x+{2j-1\over 2}\right), &m\geq 0\\
1, &m=-1\\
\prod_{j=1}^{-m-1}\left( x-{2j+1\over 2}\right)^{-1}, &m< -1.
\end{array}
\right.
$$
We obtain the well known realization of the Virasoro algebra 
in the theory of KdV hierarchy and in 2D
topological gravity \cite{DVV1,FKN}.

\medskip
\noindent
{\bf Example 2.2.} Let $V=H^*(X; {\bf C})$ where $X$ is a smooth
projective variety of complex dimension $d$ having $H^{\rm odd}(X)=0$.
The bilinear form $<~,~>$ is specified by Poincar\'e pairing
$$
<\omega_1, \omega_2> =\int_X \omega_1\wedge\omega_2.
$$
The matrix $\mu$ is diagonal in the homogeneous basis
$$
e_\alpha\in H^{2 q_\alpha}(X), ~~\mu(e_\alpha) = \left( q_\alpha -{d\over
2}\right) e_\alpha.
$$
The matrix $R=R_1$ (i.e., $R_2=R_3=\dots =0$) is the matrix of
multiplication by the first Chern class $c_1(X)$. In this case the
operators $L_m$ coincide with the Virasoro operators of \cite{Egu2}
(introduced in this paper for the particular case of quantum
cohomology of $X={\bf CP}^d$).
Recall that the chiral primary fields of the quantum cohomology of $X$
(of the topological sigma-model coupled to gravity with the target space
$X$, using physical language)
are in one-to-one correspondence with a homogeneous basis in
$H^*(X)$,
$$
\phi_1=1\in H^0(X), \quad \phi_\alpha\in H^{2 q_\alpha}(X).
$$
The topological correlators are defined by
$$
\left< \tau_{n_1} (\phi_{\alpha_1} ) \dots \tau_{n_k} (\phi_{\alpha_k})
\right>_g
$$
$$
:= \sum_{\beta\in H_2(X,{\bf Z})} q^\beta \int_{\left[\bar{\cal M}_{g,k}
(X,\beta) \right]^{\rm virt}}
c_1(L_1)^{n_1} \cup {\rm ev}_1^* (\phi_{\alpha_1}) \cup \dots \cup
c_1(L_k)^{n_k} \cup {\rm ev}_k^* (\phi_{\alpha_k})
$$
In this formula $q^\beta$ belongs to the Novikov ring ${\bf Z}$-spanned by
monomials
$q^\beta:=
q_1^{b_1} \dots q_r^{b_r}$ for a basis $q_1$, \dots, $q_r$ of $H_2(X,{\bf
Z})$ where $\beta = \sum b_i q_i$. The integration is taken over the
virtual moduli space ${\left[\bar{\cal M}_{g,k}
(X,\beta) \right]^{\rm virt}}$ of the degree $\beta\in H_2(X,{\bf Z})$ 
stable maps from
$k$-marked curves of genus $g$ to $X$; $c_1(L_i)$ is the first Chern class
of the tautological line bundle $L_i$ over ${\left[\bar{\cal M}_{g,k}
(X,\beta) \right]^{\rm virt}}$ whose fiber over each stable map is defined
by the cotangent space of underlying curve at the $i$-th marked point.
Finally, ${\rm ev}_i$ is the evalution map from ${\left[\bar{\cal M}_{g,k}
(X,\beta) \right]^{\rm virt}}$ to $X$ defined by evaluating each stable
map at the $i$-th marked point.
According to Remark 2.3 the operators $L_m$ with $m\geq -1$
are always well-defined. The full Virasoro algebra exists for even complex
dimension $d$. See also \cite{GP} where certain higher genus Virasoro
constraints were studied.

\smallskip
We define now another embedding of Virasoro algebra into the universal
enveloping of the Heisenberg algebra (\ref{jw16}). 
Let us choose an eigenvector
$e_1$ of the operator $\mu$,
$$
\mu(e_1) =\mu_1 e_1.
$$
\smallskip
\begin{lemma} The shift $a_k \mapsto \tilde a_k$ where
\beq\label{jw8}
\tilde a_k =\left\{
\begin{array}{ll}
a_k, &k\neq -2\\
a_k - e_1 {\bf 1}, &k=-2
\end{array}
\right.
\eeq
in the formula (\ref{jw3}) preserves the Virasoro commutation 
relations (\ref{su3}).
\end{lemma}

The proof is straightforward.

\medskip
The resulting Virasoro algebra will be denoted $Vir\left( V,<~,~>, \mu,
e_1,R\right)$. It depends only on the equivalence class of the data
$\left( V,<~,~>, \mu,e_1, R\right)$ under the transformations (\ref{jw4-1})
such that the operator $G_0$ respects the marked eigenvectors
$G_0(e_1^{(1)})=e_1^{(2)}$. 

\smallskip
At the end of this section we construct a free field representation of the
Virasoro algebra. This construction is analogous to \cite{Egu2}.

Let us introduce the operator-valued vector function of $z$
\beq
\phi(z;\mu) =-\sum_{k\in \bf Z} z^{-k-{1\over 2} +\mu}
{\Gamma\left( {1\over 2} -\mu + k\right)\over
\Gamma\left({1\over 2}-\mu\right)}a_k
\eeq
and the current
\beq
j(z;\mu)=\partial_z \phi(z;\mu).
\eeq
Then the Virasoro operators are the coefficients of the expansion
of the stress-tensor
\eqa
&&T(z)=\sum_{k\in {\bf Z}}L_k z^{-k-2}
={1\over 4 z^2} {\rm tr}\,\left( {1\over 4} -\mu^2\right)
\\
&&
+{1\over 2} :\left< \Gamma\left({1\over 2} -\mu +x\right) j(z; \mu-x),
e^{R\,{\buildrel\leftrightarrow \over\partial}_x} \Gamma^{-1}\left({1\over
2}+\mu+x\right) j(z;
\mu+x)\right>_{x=0}: .\nn
\eeqa
Here the operator $e^{R\, {\buildrel\leftrightarrow\over\partial}_x}$ is
defined by
the equation
\beq
\left< f(x) , e^{R\, {\buildrel\leftrightarrow\over\partial}_x}
g(x)\right>=
\left<e^{R^* \partial_x} f(x), e^{R\, \partial_x}g(x)\right> 
\eeq
where $R^*$ is the conjugate operator $<R^*(a),b>=<a,R(b)>$.
Observe that
\beq
R^* =\left(R_1+R_2+R_3+\dots\right)^* = R_1 -R_2 +R_3 -\dots.
\eeq

\setcounter{equation}{0}
\vskip 0.6cm
\section{From Frobenius manifolds to partition functions}\par
\medskip
{\it WDVV equations of associativity} is the problem of
finding of  a
quasihomogeneous, up to at most quadratic polynomial,  function $F(t)$ of
the variables $t=(t^1, \dots, t^n)$
and of a constant nondegenerate symmetric matrix
$\left(\eta^{\alpha\beta}\right)$ such that the following combinations
of the third derivatives
\beq
c_{\alpha\beta}^\gamma(t):=\eta^{\gamma\epsilon} 
\partial_\epsilon \partial_\alpha \partial_\beta F(t)
\eeq
for any $t$ are structure constants of an asociative algebra
$$
A_t ={\rm span}\, (e_1, \dots, e_n) , ~e_\alpha \cdot e_\beta =
c_{\alpha\beta}^\gamma(t) e_\gamma, ~~\alpha, \, \beta =1, \dots, n
$$
with the unity $e=e_1$ (summation w.r.t. repeated indices will be
assumed). 
Frobenius manifolds were introduced in \cite{D5} as the coordinate-free form
of WDVV. We recall now the definition.

A {\it Frobenius algebra} (over the field $k={\bf C}$)
is a pair $(A,<~,~>)$, where $A$ is a commutative
associative $k$-algebra with a unity $e$, $<~,~>$ is a symmetric
nondegenerate
{\it invariant} bilinear form $A\times A\to k$, i.e. $<a\cdot
b,c>=<a,b\cdot c>$ for any $a,\, b, \, c \in A$. A {\it gradation of
charge
$d$} on $A$ is a $k$-derivation $Q:A\to A$ such that
\beq
<Q(a),b>+<a,Q(b)>=d<a,b>, ~d\in k.
\eeq
More generally, graded of charge
$d\in k$ Frobenius algebra $(A,<~,~>)$ over a graded commutative
associative
$k$-algebra $R$ by definition is endowed with two $k$-derivations
$Q_R:R\to R$ and $Q_A:A\to A$ satisfying the properties
\eqa
&&Q_A(\alpha a) =Q_R(\alpha) a+ \alpha Q_A(a), ~\alpha\in R, ~a\in A
\\
&&
<Q_A(a), b> + <a,Q_A(b)>-Q_R<a,b> = d <a,b> , ~a,\,
b\in A.
\eeqa

A {\it Frobenius structure} of charge $d$ on the manifold
$M$ is the
structure of a Frobenius algebra on the tangent spaces $T_tM
=(A_t,<~,~>_t)$ depending (smoothly, analytically etc.) on the point $t\in
M$. It must satisfy the following axioms:
\vskip 0.3cm
\noindent
{\bf FM1.} The metric $<~,~>_t$ on $M$ is flat (but not necessarily
positive definite). Denote $\nabla$ the Levi-Civita connection for the
metric. The unity vector field $e$ must be covariantly constant, $\nabla
e=0$.

\vskip 0.3cm
\noindent
{\bf FM2.} Let $c$ be the 3-tensor $c(u,v,w):=<u\cdot v, w>$, $u,\, v,\,
w\in T_tM$. The 4-tensor $(\nabla_z c)(u,v,w)$ must be symmetric in
$u,\, v,\,  
w, \, z \in T_tM$. 

\vskip 0.3cm
\noindent
{\bf FM3.} A linear vector field $E\in Vect(M)$ must be fixed on $M$,
i.e. $\nabla\nabla E=0$, such that the derivations
$Q_{Func(M)}:=E, ~~Q_{Vect(M)}:={\rm id}+{\rm ad}_E
$
introduce in $Vect(M)$ the structure of graded Frobenius algebra of the
given charge $d$ over the graded ring $Func(M)$ of (smooth, analytic etc.)
functions on $M$. We call $E$ the {\it Euler vector field}.

\vskip 0.3cm
A manifold $M$ equipped with a Frobenius structure is called a Frobenius
manifold.
\vskip 0.3cm
Locally, in the flat coordinates $t^1, \dots, t^n$ for the metric
$<~,~>_t$, a FM is described by a solution $F(t)$ of the WDVV associativity 
equations,
where $\partial_\alpha\partial_\beta\partial_\gamma
F(t)=<\partial_\alpha\cdot \partial_\beta, \partial_\gamma>$, and vice
versa. We will call $F(t)$ {\it the potential} of the FM
(physicists call it {\it primary free energy}; in the setting of quantum 
cohomology it is called the Gromov - Witten potential \cite{KM}). 
Observe that the
unity $e$ is an eigenvector of $\nabla E$ with the eigenvalue $1$.
The
eigenfunctions of the operator $\partial_E$, $\partial_E f(t) = p\, f(t)$
are called {\it quasihomogeneous} functions of degree $p$ on the
Frobenius manifold. 
In main
applications the $(1,1)$-tensor $\nabla E$ is diagonalizable. The
potential $F(t)$ is a quasihomogeneous function of the degree $3-d$ well-defined
up to
a quadratic polynomial.

A generalization of the above definition to Frobenius supermanifolds was
proposed in \cite{KM}. We will not consider it in the present article.

An important geometrical structure on a Frobenius manifold is the
{\it deformed flat connection} $\tilde \nabla$ introduced in \cite{D3}. 
It is
defined by
the formula
\beq
\tilde \nabla_u v:= \nabla_u v + z\, u\cdot v.
\eeq
Here $u$, $v$ are two vector fields on $M$, $z$ is the parameter of the
deformation. (In \cite{Gi1} another normalization is used 
$\tilde\nabla \mapsto \hbar \tilde\nabla$, $\hbar=z^{-1}$.) We extend this
to a meromorphic connection on the direct product $M\times {\bf C}$, $z\in
{\bf C}$, by the
formula
\beq
\tilde\nabla_{d/dz} v=\partial_z v +E\cdot v -{1\over z}\mu\, v 
\eeq
with
\beq\label{zh4}
\mu:= {1\over 2}(2-d)\cdot {\bf 1} -\nabla E,
\eeq
other covariant derivatives are trivial.
Here $u$, $v$ are tangent vector fields on $M\times {\bf C}$ having zero
components along ${\bf C}\ni z$.
The curvature of $\tilde\nabla$ is equal
to zero \cite{D3}.
So, there locally exist $n$ independent
functions
$\tilde t_1(t;z), \dots, \tilde t_n(t;z)$, $z\neq 0$, such that
\beq\label{zh2}
\tilde\nabla\, d\tilde t_\alpha(t;z)=0, ~\alpha=1, \dots, n.
\eeq
We call these functions {\it deformed flat coordinates}. 

\smallskip 
\noindent{\bf Remark 3.1.} In the setting of quantum cohomology the part
of the
system (\ref{zh2}) 
not containing $\tilde\nabla_{d/dz}$ is called {\it quantum
differential equations} \cite{Gi1}. 

The component $\tilde\nabla_{d/dz}$ of the system (\ref{zh2}) for the vectors
$\xi=\nabla \tilde t$ reads
\beq\label{zh3}
z\partial_z\xi =(\mu+z {\cal U})\xi
\eeq
where
\beq\label{ad8}
{\cal U} =\left({\cal U}_\beta^\alpha(t)\right), ~~{\cal
U}_\beta^\alpha(t) :=E^\epsilon(t) c_{\epsilon\beta}^\alpha (t)
\eeq
is the matrix of the operator of multiplication by the Euler vector field.
The matrix-valued function
\beq
A(z) =\mu+z {\cal U}
\eeq
satisfies the symmetry (\ref{jw5})
\beq
\left< A(-z) a, b\right> +\left< a, A(z)b\right> =0
\eeq
for any two vectors $a$, $b$. The  equivalence class of the system 
(\ref{zh3})
w.r.t. the gauge transformations $\xi\mapsto \Theta(z)\xi$ with
$\Theta(z)$ of the form (\ref{jw6}) 
is called {\it the monodromy data at } $z=0$
of the Frobenius manifold. As was shown in \cite{D4} (for the case of
diagonalizable matrix $\mu$, the general case is similar), the monodromy
data are uniquely determined by a quadruple $\left( V, <~,~>, \mu,
R\right)$ and by a marked eigenvector $e_1:=e$ of $\mu$, $\mu(e_1)
=-{d\over 2}e_1$,  of the form described
in Section 1 above. Here $V=T_tM$ with the
bilinear form $<~,~>$ and with the antisymmetric operator $\mu$ given by
(\ref{zh4}). 
The $\mu$-nilpotent matrix $R$ gives an additional set of parameters
in the case of resonant $\mu$. The monodromy data do not depend on $t\in
M$ \cite{D4}. 

Particularly, in the setting of quantum cohomology the monodromy data
has the form of Example 2.2 above (see \cite{D4}).

From the structure (\ref{jw7}) of the fundamental system
of solutions we obtain the existence of a system of deformed flat coordinates
of the form
\beq\label{zh5}
\left( \tilde t_1(t;z), \dots, \tilde t_n(t;z)\right) =
\left( \theta_1 (t;z), \dots, \theta_n (t;z)\right) z^\mu z^R
\eeq
where the functions  $\theta_1 (t;z)$, \dots, $\theta_n (t;z)$
are analytic at $z=0$ and they satisfy the conditions
\beq
\theta_\alpha(t;0) = \eta_{\alpha\beta}t^\beta, ~~\alpha=1, \dots, n
\eeq
\beq\label{zh6}
\left< \nabla \theta_\alpha (t;-z), \nabla \theta_\beta (t; z)\right> =
\eta_{\alpha\beta}.
\eeq
In the notations of (\ref{zh5}) 
the columns of the matrix $\Theta(t;z)= \left(
\Theta_\beta^\alpha (t;z)\right)$ are the gradients
$\nabla\theta_1 (t;z)$, \dots, $\nabla\theta_n (t;z)$,
\beq\label{zh1}
\Theta_\beta^\alpha(t;z) =\nabla^\alpha 
\theta_\beta(t;z)=\delta_\beta^\alpha+O(z). 
\eeq
The deformed flat coordinates (\ref{zh5}) are determined uniquely up to
transformations of the form (\ref{jw4-1}).

Let us denote $\theta_{\alpha,p}(t)$  the coefficients of Taylor
expansions
\beq
\theta_\alpha(t;z) =\sum_{p=0}^\infty \theta_{\alpha,p}(t)z^p.
\eeq
We introduce also, due to (\ref{zh6}), the matrix-valued function
\beq
\Omega(t; z,w) =\left( \Omega_{\alpha\beta}(t;z,w)\right)
={1\over z+w} \left[ \Theta^T(t; z)\,\eta\, \Theta(t; w)-\eta\right]
\eeq
and its Taylor coefficients
\beq
\Omega_{\alpha\beta}(t;z,w) =\sum_{p,\, q=0}^\infty
\Omega_{\alpha,p; \, \beta,q}(t) z^p w^q.
\eeq
The coefficients $\Omega_{\alpha,p;\, \beta,q}(t)$ can be expressed via
scalar products of the gradients of the functions $\theta_{\alpha,m}(t)$,
$\theta_{\beta,m}(t)$ with various $m$
\beq
\Omega_{\alpha,p; \, \beta,q}(t) =\sum_{m=0} ^q (-1)^m \left< \nabla
\theta_{\alpha, p+m+1}, \nabla\theta_{\beta,q-m}\right> .
\eeq
They also satisfy the following identities
\beq\label{zh7}
\Omega_{\al,p+1;\beta,q}+\Omega_{\al,p;\beta,q+1}=
\Omega_{\al,p;\sigma,0}\,\eta^{\sigma\rho}\,\Omega_{\rho,0;\beta,q}.
\eeq
For the gradients of the entries of the matrix $\Omega(t;z,w)$ one has
\beq\label{gradomega}
\nabla\Omega_{\alpha\beta}(t;z,w) =\nabla\theta_\alpha(t;z)\cdot 
\nabla\theta_\beta(t;w).
\eeq
From this and from (\ref{zh1}) it readily follows that the functions
\beq\label{flat}
t^\alpha:=\eta^{\alpha\beta}\Omega_{\beta,0;\, 1,0}
\eeq
can serve as the flat coordinates on $M$ and
\beq
\partial_\beta t^\alpha = \delta_\beta^\alpha.
\eeq
Particularly, the unity vector field is
$$
e=e_1 ={\partial\over \partial t^1}.
$$
Moreover, if we define the function $F$ on $M$ by 
\beq\label{eff}
F(t) ={1\over 2} \left[ \Omega_{1,1; \, 1,1}(t) - 2 t^\alpha 
\Omega_{\alpha,0;\, 1,1}(t) + t^\alpha t^\beta \Omega_{\alpha,0; \, 
\beta,0}(t)\right]
\eeq
with the flat coordinates of the form (\ref{flat}) then
$$
\left< \partial_\alpha \cdot \partial_\beta, \partial_\gamma\right>
=\partial_\alpha\partial_\beta\partial_\gamma F(t).
$$
So $F(t)$ is the potential of the Frobenius manifold. It will be also
useful to know that the derivatives of the potential defined by (\ref{eff})
are
\beq\label{zh23}
\partial_\alpha F= t^\beta \Omega_{\alpha,0; \beta,0}(t) 
-\Omega_{\alpha,0;1,1}(t)
\eeq
\beq\label{zh23_1}
\partial_\alpha\partial_\beta F =\Omega_{\alpha,0; \beta,0}(t).
\eeq

Now we are able to construct the genus zero approximation to the
integrable hierarchy (\ref{su6}) and to compute the genus zero free energy
${\cal F}_0(T)$. 

The genus zero hierarchy reads
\beq\label{zh8}
\partial_{T^{\alpha,p}} t = \partial_X K^{(0)}_{\alpha,p}(t), ~~
 K^{(0)}_{\alpha,p}(t)=\nabla \theta_{\alpha,p+1}(t).
\eeq
This is a Hamiltonian hierarchy on the loop space $L(M)$ equipped with the
Poisson bracket $\{ ~,~\}_1^{(0)}$
\beq\label{zh10}
\left\{ t^\alpha(X), t^\beta(Y)\right\}_1^{(0)}=
\eta^{\alpha\beta}\delta' (X-Y).
\eeq
Here $\delta'(X)$ is the derivative of the delta-function on the circle.
The Hamiltonian reads
\beq
H_{\alpha,p}^{(0)} ={1\over 2\pi} \int_0^{2\pi}
\theta_{\alpha,p+1}\left(t(X)\right)\,dX.
\eeq
Particularly,
$$
\partial_{T^{1,0}}t=\partial_X t.
$$
So we will identify $X$ and $T^{1,0}$.

The second Hamiltonian structure of the hierarchy (\ref{zh8}), 
found in \cite{D2}, has
the form
\beq\label{zh11}
\left\{ t^\alpha(X), t^\beta(Y)\right\}_2^{(0)}
= g^{\alpha\beta}(t(X)) \delta'(X-Y) +\Gamma_\gamma^{\alpha\beta}(t)
t^\gamma_X \delta(X-Y).
\eeq
Here
\beq\label{zh9}
g^{\alpha\beta}(t) =\eta^{\alpha\gamma} {\cal U}_\gamma^\beta(t)
\eeq
is {\it the intersection form} of the Frobenius manifold $M$ 
(see \cite{D3}).
It defines on an open subset of $M$, where $\det \left(
g^{\alpha\beta}(t)\right)\neq 0$, a new flat metric. The coefficients
\beq
\Gamma_\gamma^{\alpha\beta}(t)=c_\gamma^{\alpha\epsilon}(t) \left(
{1\over 2} -\mu\right)_\epsilon^\beta
\eeq
are related to the Christoffel coefficients
$\Gamma_{\epsilon\gamma}^\beta(t)$ of the Levi-Civita connection of
the metric (\ref{zh9}) 
\beq
\Gamma_\gamma^{\alpha\beta}(t)
=-g^{\alpha\epsilon}(t)\Gamma_{\epsilon\gamma}^\beta(t).
\eeq
The Poisson brackets (\ref{zh10}) and (\ref{zh11}) 
are {\it compatible}. This means
\cite{Magri} that
any linear combination
$$
\{ ~,~\}_2^{(0)} -\lambda \{ ~,~\}_1^{(0)}
$$
with an arbitrary value of $\lambda$ is again a Poisson bracket. All these
statements follow from the general theory of \cite{DN} of Poisson brackets
of hydrodynamic type (i.e., those of the form (\ref{zh11}) 
and from the theorem
\cite{D3} 
that the metrics $g^{\alpha\beta}(t)$ and $\eta^{\alpha\beta}$ form 
{\it a flat pencil} (see also \cite{D7}). 
Particularly, denoting $\hat\nabla$
the Levi-Civita covariant gradient for the metric (\ref{zh9}), we obtain the
following important formula
\beq
\hat\nabla d\tilde t \equiv \left( \hat\nabla^\alpha \partial_\beta \tilde
t\right) =d \left( \partial_z -{1\over 2\, z}\right) \nabla \tilde t
\eeq
for an arbitrary deformed flat coordinate $\tilde t = \tilde t(t;z)$ of
the deformed connection $\tilde \nabla$.

The bihamiltonian structure (\ref{zh10}) and 
(\ref{zh11}) gives the possibility \cite{Magri}
to construct the genus zero hierarchy applying the recursion operator
\beq\label{zh20}
{\cal R} =\left( {\cal R}_\beta^\alpha\right), ~~{\cal R}_\beta^\alpha
={\cal U}_\beta^\alpha (t) +\Gamma_\gamma^{\alpha\epsilon}(t) t^\gamma_X
\eta_{\epsilon\beta}\partial_X^{-1}.
\eeq

We now specify a particular solution of the hierarchy. Due to the obvious
invariance of the equations of the hierarchy w.r.t. the transformations
$$
T^{\alpha,p}\mapsto c\, T^{\alpha,p}, ~~t\mapsto t
$$
and w.r.t. shifts along the times $T^{\beta,q}$, we specify the
symmetric solution
$t^{(0)}(T)$ as follows
\beq
\sum_{\alpha,p} \tilde T^{\alpha,p} \partial_{T^{\alpha,p}} t^{(0)}=0
\eeq
where the shifted times are defined in (\ref{su7}). 
This solution can be found
from the fixed point equation
\beq\label{zh24}
t=\nabla\Phi_T(t)
\eeq
where
\beq
\Phi_T(t)=\sum T^{\alpha,p}\theta_{\alpha,p}(t)
\eeq
in the form 
$$
t^{(0)}(T) =T^0 + \sum_{q>0} T^{\beta,q}\nabla\theta_{\beta,q}(T^0)
+\sum_{p,\, q>0} T^{\beta,q}
T^{\gamma,p}\nabla\theta_{\beta,q-1}(T^0)\cdot \nabla
\theta_{\gamma,p}(T^0)+\dots .
$$
This is a power series in $T^{\alpha,p>0}$ with coefficients 
that  depend on $T^{\alpha,0}$ through certain functions defined
on the Frobenius manifold, where we are to substitute $t^\alpha\mapsto
T^{\alpha,0}$ for any $\alpha$.

Finally, the genus zero free energy has the form
\beq\label{zh12}
{\cal F}_0(T) ={1\over 2} \sum \Omega_{\alpha,p;\, \beta,q}\left(
t^{(0)}(T)\right) \tilde T^{\alpha,p}\tilde T^{\beta,q}.
\eeq
Particularly, from this formula it follows that
\beq
\Omega_{\alpha,p; \, \gamma,r}(t)|_{t=t^{(0)}(T)}=\left<
\tau_p(\phi_\alpha)\tau_r(\phi_\gamma)  
e^{\sum \tau_q(\phi_\beta)T^{\beta,q}}\right>_0
=\partial_{T^{\alpha,p}}\partial_{T^{\gamma,r}} {\cal F}_0(T).
\label{zh22-1}
\eeq

\medskip
The genus 1 free energy has the following structure, according to 
\cite{DW}
\beq\label{zh19}
{\cal F}_1(T) =
\left[ G(t) +{1\over 24} \log\det M^\alpha_\beta(t,t_X)\right]_
{t=t^{(0)}(T), ~t_X =\partial_{T^{1,0}}t^{(0)}(T)}.
\eeq
Here the matrix $M^\alpha_\beta(t,t_X)$ is given by the formula
$$
M^\alpha_\beta(t,t_X) =c^\alpha_{\beta\gamma}(t)
t^\gamma_X
$$
and $G(t)$ is some function on the Frobenius manifold $M$. In the case
of quantum cohomology this function $G(t)$
was identified in \cite{Ge} as the
generating function of elliptic Gromov - Witten invariants. As it was
shown in \cite{DZ2}, 
this function can be computed in terms of the Frobenius
structure by a universal formula provided that 
the semisimplicity condition
for $M$ holds true.

To explain here this formula we recall first necessary constructions of
the theory \cite{D3} of
semisimple FMs.

A point $t\in M$ is called {\it
semisimple} if the
algebra on $T_tM$ is semisimple. A connected FM $M$ is
called
semisimple if it has at least one semisimple point.
Classification of semisimple FMs can be reduced, by a
nonlinear change of coordinates, to a system of ordinary differential
equations. First we will describe these new coordinates.
 
Denote $u_1(t)$, \dots, $u_n(t)$ 
the roots of the characteristic
polynomial of the operator ${\cal U}(t)$ of multiplication by the Euler
vector field
$E(t)$ ($n={\rm dim}\, M$). Denote $M^0\subset M$ the open subset where
all the roots are pairwise distinct.
It turns out \cite{D1} that $M^0$ is nonempty and
the functions $u_1(t)$, \dots, $u_n(t)$
are
independent local coordinates on $M^0$. In these coordinates
\beq\label{zh21}
\partial_i\cdot \partial_j =\delta_{ij} \partial_i, ~~{\rm
where}~\partial_i:=\partial/\partial u_i,
\eeq
and 
\beq\label{zh22}
E=\sum_i u_i \partial_i.
\eeq
The local coordinates $u_1$, \dots, $u_n$ on $M^0$ are called {\it
canonical}. 

We will now rewrite WDVV in the canonical coordinates.
Staying in a small ball on $M^0$, let us order
the canonical coordinates and choose the signs of the square roots
\beq
\psi_{i1}:=\sqrt{<\partial_i, \partial_i>}, ~~i=1, \dots, n.
\eeq
The orthonormal frame of the normalized idempotents $\partial_i$
establishes a local trivialization of the tangent bundle $TM^0$. The
deformed flat connection $\tilde \nabla$ in $TM^0$ is recasted into the
following flat connection in the trivial bundle $M^0\times {\bf C}\times
{\bf C}^n$
\beq\label{zh14}
\tilde\nabla_i\mapsto {\cal D}_i = \partial_i -z\, E_i -V_i, ~~~~
\tilde\nabla_{d/dz}\mapsto {\cal D}_z = \partial_z -U - {1\over z}V,
\eeq
other components are obvious. Here the  $n\times n$ matrices $E_i$, $U$,
$V=(V_{ij})$ read
\beq\label{zh25}
(E_i)_{kl}=\delta_{ik}\delta_{il}, ~~U={\rm diag}\, (u_1, \dots, u_n),
~V=\Psi \, \mu\, \Psi^{-1}=-V^T
\eeq
where the matrix $\Psi=(\psi_{i\alpha})$ satisfying $\Psi^T\Psi=\eta$ is
defined by
\beq
\psi_{i\alpha} := 
\psi_{i1}^{-1}\,
\frac{\partial t_\alpha}{\partial u_i},
~~i, \, \alpha=1, \dots, n.
\eeq
The skew-symmetric matrices $V_i$ are determined by the equations
\beq\label{zh26}
[U,V_i]=[E_i,V].
\eeq

Flatness of the connection (\ref{zh14}) 
$$
[{\cal D}_i,{\cal D}_j]=0, ~~[{\cal D}_i,{\cal D}_z]=0
$$
reads as the system of commuting
time-dependent Hamiltonian flows on the Lie algebra $so(n)\ni V$ equipped
with the standard linear Poisson bracket
\beq\label{zh13}
\frac{\pal V}{\pal u_i}=\left\{ V,H_i(V;u)\right\} ,~~i=1, \dots, n
\eeq
with the quadratic Hamiltonians
\beq\label{zh15}
H_i(V;u) ={1\over 2} \sum_{j\neq i} {V_{ij}^2\over u_i-u_j}, ~i=1, \dots,
n.
\eeq
The monodromy of the operator $\tilde\nabla_{d/dz}$ (i.e., the monodromy
at the origin, the Stokes matrix, and the central connection matrix, see
definitions in \cite{D3,D4}) does not change with small variations of a
point $u=
(u_1, \dots, u_n)\in M$. Using this isomonodromicity property one
can parametrize semisimple Frobenius manifolds by the monodromy data
of the deformed flat connection (see \cite{D4}). 
Recall that in our construction
of the Virasoro constraints we use just the  monodromy at the origin
defined
for an arbitrary Frobenius manifold, not only for a semisimple one.

We are ready now to give the formula for the function $G(t)$ called in
\cite{DZ2} the
$G$-function of the Frobenius manifold. We define, following \cite{JM}, 
the
tau-function of the isomonodromy deformation (\ref{zh13}) by a quadrature
\beq\label{zh16}
d\log \tau = \sum_{i=1}^n H_i(V(u);u) du_i.
\eeq
The $G$-function of the FM is defined by
\beq\label{zh17}
G=\log ({\tau/ J^{1/24}})
\eeq
where 
\beq\label{zh18}
J=\det \left( \partial t^\alpha /\partial u_i\right) =\pm \prod_{i=1}^n
\psi_{i1}(u).
\eeq
This function is defined on some covering $\hat M^0$ of $M^0\subset M$.
As  was proved in \cite{DZ2}, if the Frobenius manifold $M$ is the
quantum cohomology of a smooth projective variety $X$, and $M$
is semisimple, then the $G$-function of $M$ coincides with the generating
function of elliptic Gromov - Witten invariants.

The genus 0 hierarchy (\ref{zh8}) 
can be extended to the genus 1 (see (\ref{su15})). The
explicit formula for the correction $K^{(1)}(t,t_X, t_{XX})$ together with
an extension of the bihamiltonian structure (\ref{zh10}), 
(\ref{zh11}) was obtained in \cite{DZ2}.
We will not give here these formulae.

We are now ready to return to Virasoro constraints. Let us consider the
algebra
$Vir\left( V,<~,~>, \mu,
e_1,R\right)$ where $\left( V,<~,~>, \mu,
e_1,R\right)$ is the monodromy at the origin of the Frobenius manifold
$M$. The realization (\ref{jw3}), (\ref{jw9}) 
 of the algebra, where one is to do the shift (\ref{jw8})
corresponding to (\ref{su7}), acts on the algebras of ``functions''
$$
{\bf A}_0 := Func(M)\otimes {\bf C}[[T^{\alpha,p>0}]]
$$
or on
$$
{\bf A}_1 := Func(\hat M^0)\otimes {\bf C}[[T^{\alpha,p>0}]].
$$
Here we embed the algebras of smooth functions $Func(M)$ or $Func(\hat
M^0)$ into ${\bf A}_0$ or ${\bf A}_1$ resp. mapping any function $f(t)$ 
to $f(T^0)$. Observe that ${\cal F}_0(T)\in {\bf A}_0$,  ${\cal F}_1(T)\in
{\bf A}_1$.

So, our Main Theorem says that , for an arbitrary Frobenius manifold
the function ${\cal F}_0(T)$ given by (\ref{zh12})  
satisfies the differential
equations (\ref{su16}),
and for an arbitrary semisimple Frobenius manifold the function ${\cal
F}_1(T)$ given by (\ref{zh19}) satisfies the differential equations 
(\ref{su17}).

The proof of this theorem will be given in the next section.

\setcounter{equation}{0}
\vskip 0.5cm
\section{From symmetries of the hierarchy to Virasoro
constraints}

\noindent
{\bf Definition 4.1.} A {\it symmetry} of the hierarchy (\ref{zh8}) is an
evolutionary system 
$$
\partial_s t =S\left( t,t_X, \dots; X, T\right)
$$
commuting with all equations of the hierarchy
$$
\partial_s \partial_{T^{\alpha,p}}=\partial_{T^{\alpha,p}}\partial_s.
$$
All the symmetries of the hierarchy form a Lie algebra w.r.t. the
commutator.

\smallskip
\noindent
{\bf Example 4.1.} Any of the flows of the hierarchy is a symmetry.

\smallskip
\noindent
{\bf Example 4.2.} The rescaling combined with the shift along $T^{1,1}$ 
$$
\partial_s t =\sum \tilde T^{\alpha,p} \partial_{T^{\alpha,p}} t
$$
is a symmetry of the hierarchy.

Recall that the solution $t^{(0)}(T)$ is invariant w.r.t. this symmetry.
From this it readily follows that the genus zero free energy is
a homogeneous function of  degree 2 of the variables $\tilde
T^{\alpha,p}$. Indeed, the coefficients of the quadratic form (\ref{zh12}) 
being functions on $t^{(0)}(T)$ are
constant along the vector field $\sum {\tilde T}^{\al,p}\,\pal_{T^{\al,p}}$.

\smallskip
\noindent
{\bf Example 4.3.} For {\it Galilean symmetry} 
\beq\label{bo1}
S_{-1} = \sum_{p\geq 1} \tilde T^{\alpha,p}\partial_{T^{\alpha,p-1}}t +
e_1.
\eeq
Here, as above, $e_1$ is the unity vector field on the Frobenius manifold.

\medskip
To produce other symmetries we will apply the recursion relation (\ref{zh20})
 to
the symmetry (\ref{bo1}). We obtain a new symmetry
\eqa
&&S_0 ={\cal R}\, S_{-1} = \left[
\sum_{p\geq 0}\left< \left(p+\mu+{1\over 2}\right) \tilde T^p,
\partial_{T^p}\right>
+\sum_{p\geq 1}\sum_{1\leq r\leq p}
\left< R_r \tilde T^p, \partial_{T^{p-r}}\right>\right]\, t 
\nn\\
&&\quad +E(t).
\eeqa
As we will see, this symmetry corresponds to the quasihomogeneity
transformations generated by the Euler vector field $E(t)$.

Iterating this process we obtain symmetries with non-local terms in the
r.h.s.
\eqa
&&S_1 =\left[ \sum_{p\ge 0}
 \left< \left( p+\mu+{1\over 2}\right)\left( p+\mu+{3\over
2}\right) \tilde T^p, \partial_{T^{p+1}}\right> \right.
\nn\\
&&\quad\left.
+2\sum_{p\geq 0}\sum_{1\leq r\leq p+1}
\left< R_r (p+\mu+1)\tilde T^p, \partial_{T^{p-r+1}}\right>
+\sum_{p\geq 1}\sum_{2\leq r\leq p+1}
\left< R_{r,2}\tilde T^p, \partial_{T^{p-r+1}}\right>\right]\, t
\nn\\
&&\quad
+\sum \partial_{T^{\alpha,0}}t\,\left({1\over 4}-\mu^2\right)_\beta^\alpha 
\partial_X^{-1} t^\beta +E^2(t),
\\
&&S_2=\left[ \sum_{p\geq 0} \left< \left( p+\mu+{1\over 2}\right)
 \left( p+\mu+{3\over 2}\right)  \left( p+\mu+{5\over 2}\right)
\tilde T^p , \partial_{T^{p+2}}\right>
\right.
\nn\\
&&\quad
+\sum_{p\geq 0}\sum_{1\leq r\leq p+2}
\left< R_r \left[ 3  \left( p+\mu+{1\over 2}\right)^2 + 6  \left(
p+\mu+{1\over 2}\right)+2\right] \tilde T^p,
\partial_{T^{p-r+2}}\right>
\nn\\
&&\quad
+3\sum_{p\geq 0} \sum_{2\leq r\leq p+2} \left<
R_{r,2}  \left( p+\mu+{3\over 2}\right)\tilde T^p,
\partial_{T^{p-r+2}}\right>
\nn\\
&&\quad
\left.
+\sum_{p\geq 1}\sum_{3\leq r\leq p+2}
\left< R_{r,3}\tilde T^p, \partial_{T^{p-r+2}}\right> \right]\, t
\nn\\
&&\quad
+\partial_{T^{\alpha,1}} t\,\left( {3\over 8} +{1\over 4}\mu -{3\over
2}\mu^2 -\mu^3\right)_\beta^\alpha \partial_X^{-1} t^\beta
+\partial_{T^{\alpha,0}}t\,\left[R_1 \left({1\over
4}-3\mu-3\mu^2\right)\right]_\beta^\alpha 
\partial_X^{-1} t^\beta
\nn\\
&&\quad
+\partial_{T^{\alpha,0}} t \, \left( {3\over 8} -{1\over 4}\mu-{3\over
2}\mu^2 + \mu^3\right) _\beta^\alpha \eta^{\beta\gamma}\partial_X^{-1}
\partial_\gamma F(t) +E^3(t).
\eeqa

\begin{prop}\label{p4-1} The vector fields $S_m$ with $-1\leq m\leq 2$
are symmetries of the hierarchy (\ref{zh8}).
\end{prop}

\medskip
Actually, there is an unpleasant problem to choose the integration
constants in the derivation of these formulae. So we are to prove directly
that $S_m$ for $-1\leq m \leq 2$ are symmetries of the hierarchy. To do
this we will first formulate the following simple observation
\cite{KM}.
\smallskip
\begin{lemma}\label{L4-1} 
On an arbitrary Frobenius manifold the vector fields
\beq\label{bo2}
v_m =E^{m+1}(t), ~~m\geq 0, ~~v_{-1}=e
\eeq
satisfy the Virasoro commutation relations
\beq
[v_i, v_j]=(i-j) v_{i+j}, ~~i, \, j\geq -1
\eeq
with respect to the Lie bracket.
\end{lemma}

\medskip
\noindent
{\bf Remark 4.1.} On a semisimple Frobenius manifold the claim of Lemma
\ref{L4-1}
is obvious since, as it follows from (\ref{zh21}), (\ref{zh22})
$$
v_m = \sum_{i=1}^n u_i^{m+1}\partial_i.
$$
In this case the vector fields $v_m$ are defined also for any $m\in {\bf
Z}$ outside the discriminant of the Frobenius manifold (see \cite{D4}) 
and they
satisfy the commutation relations of the Virasoro algebra with zero
central charge.

\medskip
It turns that (half of) the Virasoro algebra described in
Lemma \ref{L4-1}
acts on the polynomial ring generated by the functions
$\Omega_{\alpha,p;\, \beta,q}(t)$. This follows from
\smallskip

\begin{prop}\label{p4-2} The derivatives of the coefficients
$\Omega_{\alpha,p;\, \beta,q}(t)$ along the vector fields (\ref{bo2})
are polynomials of the same functions with constant coefficients
depending only on the monodromy at the origin $\eta$, $\mu$, $R$. These
derivatives
are uniquely determined by the following formulae for the derivatives
of the generating matrix-valued function $\Omega=\Omega(t;z,w)$
\eqa
&&\pal_e\Omega=(z+w)\,\wo(z,w),\label{Wm1}\\
&&\pal_E\Omega=\left(D_z+{1\over 2}\right)^T\,\wo(z,w)+
\wo(z,w)\,\left(D_w+{1\over 2}\right),\label{W0}\\
&&\pal_{E^2}\Omega=
\left(D_z+{1\over 2}\right)^T\,\left(D_z+{3\over 2}\right)^T\,z^{-1}\,
\wo(z,w)
\nn\\
&&\quad+
\wo(z,w)\,w^{-1}\,\left(D_w+{3\over 2}\right)\,
\left(D_w+{1\over 2}\right)
\nn\\
&&\quad+\wo(z,0)\,\left({1\over 4}-\mu^2\right)
\,\eta^{-1}\,\wo(0,w),\label{W1}\\
&&\pal_{E^3}\Omega=
\left(D_z+{1\over 2}\right)^T\,\left(D_z+{3\over 2}\right)^T\,
\left(D_z+{5\over 2}\right)^T\,z^{-2}\,
\wo(z,w)\nn\\
&&\quad+
\wo(z,w)\,w^{-2}\,\left(D_w+{5\over 2}\right)\,
\left(D_w+{3\over 2}\right)\,\left(D_w+{1\over 2}\right)\nn\\
&&\quad+(\partial_w\wo)(z,0)\,\eta^{-1}\,\left({1\over 2}+\mu\right)^T\,
\left({1\over 2}-\mu\right)^T\,\left({3\over 2}-\mu\right)^T\,\wo(0,w)\nn\\
&&\quad+\wo(z,0)\,\left({3\over 2}-\mu\right)\,
\left({1\over 2}-\mu\right)\,\left({1\over 2}+\mu\right)\,\eta^{-1}\,
(\partial_z\wo)(0,w)\nn\\
&&\quad+\wo(z,0)\,\left({1\over 4}+3\mu-3\,\mu^2\right)\,R_1\,\eta^{-1}\,
\wo(0,w).\label{W2}
\eeqa
where 
\beq
\tilde\Omega(z,w):= \Omega(t;z,w) +{\eta\over z+w},
\eeq
the right operator $D_w$ acts on the matrices as
follows
\beq
\Omega D_w = w{\partial \Omega\over \partial w} + \Omega B(w),
\eeq
the operator $D_z^T$ acts by a similar formula on the left,
\beq
D_z^T\Omega =z{\partial\Omega\over\partial z} +B^T(z)\Omega
\eeq
and the matrix-valued polynomial $B$ is defined by
\beq
B(z):=\mu + \sum_{k\geq 1} R_k z^k.
\eeq
\end{prop}
\pf From the definition (\ref{zh1}) of the matrix $\Theta(t;z)$ 
and the definition (\ref{zh2}) of the deformed flat coordinates
it follows that
$$
\partial_v \Theta =z\, {\cal C}(v)\Theta
$$
where ${\cal C}(v)$ is the operator of multiplication by the vector field
$v$ tangent to $M$. From this we obtain
\beq\label{du1}
\partial_E \Theta =z {\cal U}\Theta,
\eeq
By using (\ref{zh3}) we get
\beq\label{du2}
\partial_E \Theta =\Theta D_z -\mu\, \Theta.
\eeq
Besides, from the definition it follows that
\beq
\Theta^T(t;z)\eta\Theta(t;w)-\eta = (z+w)\Omega(t;z,w).\label{du3}
\eeq
We will omit the explicit dependence on $t\in M$ in subsequent
calculations. 
 
The formula (\ref{Wm1}) follows from (\ref{du3}) and from 
$\partial_e\Theta =
z\Theta$. 
Differentiating (\ref{du3}) along $E$ and applying 
(\ref{du2}) we arrive at the
formula (\ref{W0}). We can compute the same derivative in an alternative
way
\eqa
&&\partial_E\left( \Theta^T(z)\eta\Theta(w)\right)
=z\, \Theta^T(z)\,{\cal U}^T\eta\,\Theta(w) + w\, \Theta^T(z)\,\eta\,{\cal
U}\,\Theta(w)
\nn\\
&&= (z+w)\Theta^T(z)\,\eta\,{\cal U}\,\Theta(w)=
{z+w\over w} \Theta^T(z)\,\eta\,\partial_E\Theta(w).
\eeqa
Applying (\ref{du2}) and comparing with 
$$
\partial_E\left( \Theta^T(z)\eta\Theta(w)\right)=
D_z^T\Theta(z)\,\eta\,\Theta(w)+\Theta^T(z)\,\eta\,\Theta(w)D_w
$$
we arrive at
\eqa
&&\Theta^T(z)\eta\mu\Theta(w) ={z\Theta^T(z)\eta\Theta(w)D_w
-w D_z^T \Theta^T(z)\eta\Theta(w)\over z+w}
\nn\\
&&=z\,\wo(z,w)\,D_w-w\,D_z^T\wo(z,w).\label{du4}
\eeqa
Now, to derive (\ref{W1}) we differentiate (\ref{du3}) 
along $E^2$. This can be written
as follows:
\eqa
&&\partial_{E^2} \Theta^T(z)\eta\Theta(w)=
z \Theta^T(z)\left[{\cal U}^2\right]^T\eta\Theta(w)
+w \Theta^T(z)\eta{\cal U}^2 \Theta(w) 
\nn\\
&&=(z+w) \Theta^T(z){\cal 
U}^T\eta{\cal U}\Theta(w)
={z+w\over
zw}\partial_E\Theta^T(z)\eta\partial_E\Theta(w).\nn
\eeqa
Substituting (\ref{du3}) and using (\ref{du4}) we obtain (\ref{W1}).

To prove the formula (\ref{W2}) we need the following lemma:
\begin{lemma} The following formulae hold true:
\eqa
&&\pal_{E^2}\Theta(z)=\Theta(z) (D_z-{1\over 2})(D_z-{3\over 2})\,z^{-1}  
+\cc\,({1\over 4}-\mu^2)\Theta(z)
\nn\\
&&\quad+z({1\over 2}-\mu)({3\over 2}-\mu)\,\eta^{-1}\,(\pal_w\Omega(z,0))^T
-2(\mu-1)R_1\Theta(z)
\nn\\
&&\quad-z^{-1}\,({1\over 2}-\mu)({3\over 2}-\mu),
\label{du5}
\\
&&\nn\\
&&\Theta^T(z)\,\eta\,\mu\,\cc=-z^{-1}(D_z^T-1)\, \Theta^T(z)\,\eta+
z\,\pal_w\Omega(z,0)\,(1+\mu)
\nn\\
&&\quad+\Theta^T(z)\,\eta\,R_1-z^{-1}\eta\,(1+\mu).\label{du10}
\eeqa
where the matrix $\cc$ is defined as
$$
\cc=\Omega(t;0,0)\eta^{-1}.
$$
\end{lemma}
\pf From the definition (\ref{ad8}) of  $\uu$  and (\ref{W0}) we obtain
\eqa
&&\uu=(1-\mu)\,\cc+\cc\,\mu+R_1,\label{du6}
\\
&&\pal_E{\cal U}=(1-\mu)\,{\cal U}+{\cal U}\,\mu,\label{du7}
\eeqa
By using (\ref{zh7}) we get
\beq
\cc\,\Theta=z^{-1}\,(\Theta-1)+z\,\eta^{-1}\,(\pal_w\Omega(z,0))^T.
\label{du8}
\eeq
From (\ref{du1}) and (\ref{du7}) it follows that
\eqa
&&\pal_E\pal_E\Theta=z\,\uu\,\pal_E\Theta+z\,(\pal_E\uu)\,\Theta
=z\,\pal_{E^2}\Theta+z\,((1-\mu)\,{\cal U}+{\cal U}\,\mu)\,\Theta
\nn\\
&&=z\,\pal_{E^2}\Theta+2\,(1-\mu)\,\pal_E\Theta-z\,(1-\mu)\,\uu\,\Theta
+z\,\uu\,\mu\,\Theta.
\eeqa
Insert (\ref{du6}) into the above expression we arrive at
\eqa
&&\pal_E\pal_E\Theta=z\,\pal_{E^2}\Theta+2\,(1-\mu)\,\pal_E\Theta-
z(1-\mu)[(1-\mu)\,\cc+\cc\,\mu+R_1]\,\Theta
\nn\\
&&\quad+z\,[(1-\mu)\,\cc+\cc\,\mu+R_1]\,\mu\,\Theta
\nn\\
&&=z\,\pal_{E^2}\Theta+2\,(1-\mu)\,\pal_E\Theta-z\,(1-\mu)^2\,\cc\,\Theta
+z\,\cc\,\mu^2\,\Theta+2\,z\,(\mu-1)R_1\,\Theta.\nn
\eeqa
From the above expression and by using the formula (\ref{du2}) we obtain
\eqa
&&z\pal_{E^2}\Theta=(\Theta\,D_z-\mu\,\Theta)\,D_z-\mu\,
(\Theta\,D_z-\mu\,\Theta)-2(1-\mu)(\Theta\,D_z-\mu\,\Theta)
\nn\\
&&\quad+z(1-\mu)^2\,\cc\,\Theta-
z\,\cc\,\mu^2\,\Theta-2\,z\,(\mu-1)\,R_1\,\Theta
\nn\\
&&=\Theta\,D_z^2-2\,\Theta\,D_z+(2\mu-\mu^2)\,\Theta
+z\,\cc\,({1\over 4}-\mu^2)\,\Theta
\nn\\
&&\quad+z\,({1\over 2}-
\mu)({3\over 2}-\mu)\,\cc\,\Theta-2\,z\,(\mu-1)\,R_1\,\Theta.\nn
\eeqa
By using formula (\ref{du8}) we finally get
\eqa
&&z\pal_{E^2}\Theta=
\Theta\,D_z^2-2\,\Theta\,D_z+(2\mu-\mu^2)\,\Theta
+z\,\cc\,({1\over 4}-\mu^2)\,\Theta
+({1\over 2}-
\mu)({3\over 2}-\mu)\,(\Theta-1)
\nn\\
&&\quad+z^2\,({1\over 2}-
\mu)({3\over 2}-\mu)\,\eta^{-1}\,(\pal_w\Omega(z,0))^T
-2\,z\,(\mu-1)R_1\Theta
\nn\\
&&=\Theta\,(D_z-{1\over 2})(D_z-{3\over 2})
+z\,\cc\,({1\over 4}-\mu^2)\,\Theta
\nn\\
&&\quad+z^2\,({1\over 2}-
\mu)({3\over 2}-\mu)\,\eta^{-1}\,(\pal_w\Omega(z,0))^T
-2\,z\,(\mu-1)R_1\Theta-({1\over 2}-
\mu)({3\over 2}-\mu),\nn
\eeqa
which leads to formula (\ref{du5}).

To prove formula (\ref{du10}), let's note that
\beq
\frac{\pal\Theta(z)}{\pal z}|_{z=0}=\cc,\quad
\Theta(z)|_{z=0}=I.\label{du11}
\eeq
Now differentiating (\ref{du4}) w.r.t. $w$ and putting $w=0$ we arrive at
(\ref{du10}).
\epl

Now let's prove formula (\ref{W2}). Using (\ref{du1}) and (\ref{du3})
 we get
\eqa
&&\pal_{E^3}\Omega=(z+w)^{-1}\pal_{E^3}(\Theta^T(z)\,\eta\,\Theta(w))
\nn\\
&&=(z+w)^{-1}\left(z\,\Theta^T(z)\,(\uu^3)^T\eta\,\Theta(w)+
w\,\Theta^T(z)\,\uu^3\,\eta\,\Theta(w)\right)
\nn\\
&&=(z+w)^{-1}\left[\pal_{E^2}\left(\Theta^T(z)\,\eta\,\uu\,\Theta(w)\right)-
\Theta^T(z)\,\eta\,\left(\pal_{E^2}\uu\right)\,\Theta(w)\right]
\nn\\
&&= (z+w)^{-1}\left[(z+w)^{-1}\,
\pal_{E^2}\pal_E\left(\Theta^T(z)\,\eta\,\Theta(w)\right)-
\Theta^T(z)\,\eta\,\left(\pal_{E^2}\uu\right)\,\Theta(w)\right].\nn
\eeqa
From (\ref{du3}), (\ref{du6}) and the fact
$$
\pal_{E^2}\cc=\uu^2
$$
we obtain
\eqa
&&\pal_{E^3}\Omega=(z+w)^{-1}\pal_{E^2}\pal_E\Omega-(z+w)^{-1}\Theta^T(z)
\,\eta\,\left[(1-\mu)\,\uu^2+\uu^2\,\mu\right]\,\Theta(w)
\nn\\
&&=(z+w)^{-1}\left(\pal_{E^2}\pal_E\Omega-2\,\pal_{E^2}\Omega+
w^{-1}\Theta^T(z)\,\eta\,(\frac12+\mu)\,\pal_{E^2}\Theta(w)\right.
\nn\\
&&\quad+\left.
z^{-1}\left(\pal_{E^2}\Theta^T(z)\right)\,\eta\,(\frac12-\mu)\,
\Theta(w)\right).\nn
\eeqa
By using formulae (\ref{W0}), (\ref{W1}), (\ref{du3}), (\ref{du4}),
(\ref{du5}) and (\ref{du10}) we obtain  (\ref{W2})
from the above formula by a straightforward calculation.
\epp
\medskip

\begin{cor}
(see \cite{D4}, Exercise 2.11.) In the flat coordinates (\ref{flat}) the 
Euler vector field reads
\beq\label{eulero}
E(t) =\sum\left( 1+\mu_1-\mu\right)_\beta^\alpha t^\beta \partial_\alpha
+\left( R_1\right)_1^\alpha \partial_\alpha
\eeq
(recall that $d=-2 \mu_1$). The potential $F(t)$ defined by (\ref{eff})
satisfies the following quasihomogeneity condition
\beq\label{quasiquasi}
\partial_E F(t) = (3+2\mu_1) F(t) + {1\over 2} \left< t, R_1t\right>
+\left< e_1, R_2t\right> -{1\over 2} \left< e_1, R_3e_1\right>.
\eeq
\end{cor}

\medskip
By using Proposition \ref{p4-2}, it's straightforward to prove 
Proposition \ref{p4-1}, we omit the derivation here. 
Instead,
in order to prove the genus zero Virasoro constraints,  we will
simultaneously 
show that
the solution $t^{(0)}(T)$ is  a stationary point of the symmetries 
$S_m$ for $-1\leq m \leq 2$. To be more specific, we are to say how 
to define the nonlocal terms in $S_1$ and $S_2$. Using (\ref{zh22-1})
and (\ref{zh23})
we define the nonlocal terms as follows
\eqa\label{bo3}
&&\partial_X^{-1}t^\beta 
:=\eta^{\beta\epsilon}\partial_{T^{\epsilon,0}}{\cal F}_0
\\
&&\partial_X^{-1} \partial_\gamma F := \partial_{T^{\gamma,1}} {\cal F}_0.
\label{bo4}
\eeqa

\begin{prop} For $-1\leq m\leq 2$ we have
\beq\label{bo5}
S_m \big|_{t=t^{(0)}(T)}=0
\eeq
\beq
{\cal A}_{m,0}=0
\eeq
where in $S_1$ and $S_2$ one is to substitute (\ref{bo3}) and (\ref{bo4}).
\end{prop}
\pf Let us begin with $S_{-1}$. Differentiating (\ref{zh24})
w.r.t. $X=T^{1,0}$
and using the explicit form (\ref{zh8}) 
of $\partial_{T^{\alpha,p-1}}$ we obtain
\beq
\partial_X t =\nabla \theta_{1,0} + \sum T^{\alpha,p} \partial_X
(\nabla \theta_{\alpha,p}) = e+\sum T^{\alpha,p}
\partial_{T^{\alpha,p-1}}t.
\eeq
Here we used $\nabla \theta_{1,0}=e$. We obtain (\ref{bo5}) for $m=-1$.  

Let us now derive the first Virasoro constraint (i.e., the string
equation)
\beq\label{bo6}
{\cal A}_{-1,0} \equiv\sum \tilde T^{\alpha,p+1} \partial_{T^{\alpha,p}}
{\cal F}_0(T) +{1\over 2} \left< T^0, T^0\right> =0.
\eeq
To this end we apply the operator 
$$
\sum \tilde T^{\alpha,p+1} \partial_{T^{\alpha,p}}
$$
to the function ${\cal F}_0(T)$ (see
(\ref{zh12})). Since $t^{(0)}(T)$ is a stationary point of the
symmetry $S_{-1}$,
the coefficients $\Omega_{\alpha,p;\, \beta,q}(t^{(0)}(T))$ are constants
along $S_{-1}$.  Using this and the formula (\ref{Wm1})
we easily arrive at (\ref{bo6}).

Let us proceed with the derivation of the $S_0$ symmetry of $t^{(0)}$
and of the correspondent $L_0$ Virasoro constraint. This is the first case
where we are to fix certain integration constants (see the formula 
(\ref{zh20})
for the recursion operator). 

First, applying the equation (\ref{W0}) 
to $\partial_E \Omega_{\beta,p;\, \gamma,
0}(t)$, multiplying the result by $\eta^{\alpha\gamma}$ and
differentiating along $\partial_\nu$, we obtain
\eqa
&&g^{\alpha\gamma} \partial_\gamma\partial_\nu 
\theta_{\beta,p} + \Gamma_\nu^{\alpha\gamma}\partial_\gamma
\theta_{\beta,p}
\nn\\
&&
=\left(p+\mu+{1\over 2}\right)_\beta^\gamma
\partial^\alpha\partial_\nu \theta_{\gamma,p+1} + \sum_{k=0}^{p-1}
(R_{p-k})^\gamma_\beta \,\partial^\alpha\partial_\nu
\theta_{\gamma,k+1}.
\eeqa
Multiplying this by $t_X^\nu$ we obtain, for $p>0$
\eqa
&&g^{\alpha\gamma}\partial_{T^{\beta,p-1}}t_\gamma + t^\nu_X 
\Gamma_\nu^{\alpha\gamma} \partial_\gamma\theta_{\beta,p}
\nn\\
&&
=\left( p+\mu+{1\over 2}\right)_\beta^\gamma
\partial_{T^{\gamma,p}}t^\alpha +\sum_{k=0}^{p-1} (R_{p-k})_\beta^\gamma
\,\partial_{T^{\gamma,k}}t^\alpha.
\eeqa
Multiplying by $\tilde T^{\beta,p}$ and using $S_{-1}$ we arrive, after
summation w.r.t. $\beta$ and $p\geq 1$, at
\eqa
&&-g^{\alpha\gamma}\eta_{\gamma\, 1} +t^\nu_X
\Gamma_\nu^{\alpha\gamma}
\sum_{p\geq 1} \tilde T^{\beta,p}\,\partial_\gamma \theta_{\beta,p}
\nn\\
&&
=\sum_{p\geq 1} \tilde T^{\beta,p}\left[
\left( p+\mu+{1\over 2} \right)_\beta^\gamma \partial_{T^{\gamma,k}}
+\sum_{k=0}^{p-1} 
(R_{p-k})_\beta^\gamma\, \partial_{T^{\gamma,k}}\right]\,t^\alpha.
\eeqa
From the definition of $g^{\alpha\beta}$ it follows that the first term
in the l.h.s. equals $E^\alpha(t)$. The sum in the l.h.s. is equal to
$- \tilde T^{\beta,0} \partial_\gamma \theta_{\beta,0} = -\tilde
T^{\beta,0} \eta_{\beta\gamma}$ due to the specification $\nabla
\Phi_{\tilde T}(t)=0$ of the solution $t^{(0)}$. So in the l.h.s. we
obtain
\eqa
&&-E^\alpha -\tilde T^{\beta,0}t^\epsilon_X \Gamma_\epsilon^{\alpha\gamma}
\eta_{\gamma\beta}
=- E^\alpha - \tilde T^{\beta,0} t^\nu_X c_{\nu\gamma}^\alpha
\left({1\over 2}+\mu\right) _\beta^\gamma 
\nn\\
&&= -E^\alpha -\tilde T^{\beta,0} \left( {1\over 2}+\mu\right)_\beta^\gamma
\partial_{T^{\gamma,0}}t^\alpha.
\eeqa
This proves that $S_0\big |_{t^{(0)}}=0$. Thus, applying 
the operator 
$$
\sum_{p\ge 0} \left< \left( p+\mu+{1\over 2}\right) \tilde T^p,
\partial_{T^p}\right>
+\sum_{p\geq 1} \sum_{1\leq r\leq p} \left< R_r \tilde T^p,
\partial_{T^{p-r}}\right>
$$
to the function (\ref{zh12}), 
we may assume that the coefficients
$\Omega_{\alpha,p;\, \beta,q}$ are constants along the symmetry 
$S_0$. After simple calculations with the use of (\ref{W0})
we obtain the genus zero $L_0$ Virasoro constraint
\eqa
&&\left[\sum_{p\ge 0} \left< \left( p+\mu+{1\over 2}\right) \tilde T^p,
\partial_{T^p}\right>
+\sum_{p\geq 1} \sum_{1\leq r\leq p} \left< R_r \tilde T^p,
\partial_{T^{p-r}}\right>\right]\, {\cal F}_0
\nn\\
&&
+{1\over 2} \sum_{p, \, q} (-1)^q \left< R_{p+q+1}\tilde T^p , \tilde
T^q\right> =0.
\eeqa

Note that the term $1/4 {\rm tr}\,(1/4 -  \mu^2)$ in $L_0$ (see (\ref{jw11}))
does not enter in this equation. It will appear in the genus 1 Virasoro
constraint.

Similar calculations complete the proof of Proposition.
\epp

\medskip
\begin{cor}
The genus zero Virasoro constraint 
$$
{\cal A}_{m,0}=0
$$
hold true for any $m\geq -1$.
\end{cor}

\medskip
So, the first part of Main Theorem is proved.

\smallskip
\noindent
{\bf Remark 4.2.} In \cite{Egu3} it was actualy proved  only that the
$\partial_{T^{\alpha,0}}$-derivatives of the genus zero Virasoro
constraints hold true (only the case $R=R_1$ valid in quantum cohomology
was under consideration). From this the authors of \cite{Egu3} 
infer the validity
of the Virasoro constraints choosing zero the integration constant.
A posteriori, in this case, $R=R_1$, the integration constant is zero
indeed. However, for a general Frobenius manifold the approach of \cite{Egu3}
does not work since the integration constant does not vanish.

\medskip
We proceed now to the genus one Virasoro constraints. Starting from this
point we assume the Frobenius manifold to be semisimple.
\smallskip
\begin{prop}
The derivatives of the $G$-function along the
powers
of the Euler vector field are given by the following formulae
\eqa
&&\partial_e G=0,\label{bo7}
\\
&&
\partial_E G= {n\, d\over 48} -{1\over 4} {\rm tr}\, \mu^2,
\label{bo8}
\\
&&
\partial_{E^k}G =-{1\over 4}{\rm tr}\,\left(\mu \, (\mu\, {\cal U}^{k-1}
+{\cal U}\mu\, {\cal U}^{k-2} + \dots + {\cal U}^{k-1}\mu)\right)
\nn\\
&&\quad
-{1\over 24} \left< (\mu\, {\cal U}^{k-2}+ {\cal U}\mu\, {\cal U}^{k-3}+
\dots + {\cal U}^{k-2}\mu)\,E -{d\over 2}\, {\cal U}^{k-2}\,E, H\right>,
\label{bo9}
\\
&&\quad  k\geq 2\nn
\eeqa
where
\beq\label{bo10}
H=c_\nu^{\nu\alpha}\partial_\alpha.
\eeq
\end{prop}
\pf The formula (\ref{bo7}) makes part of the definition of the
$G$-function. The formula (\ref{bo8}) was proved in \cite{DZ2}. 
To prove (\ref{bo9}) we first compute
the logarithmic derivatives of the isomonodromic tau-function along $E^k$.
Using (\ref{zh25}) and (\ref{zh16}) we have
\eqa
&&\partial_{E^k} \log \tau = {1\over 2}\sum_{i\neq j} {u_i^k\over u_i-u_j}
V_{ij}^2
={1\over 4} \sum_{i\neq j} {u_i^k-u_j^k\over u_i-u_j} V_{ij}^2
\nn\\
&&
={1\over 4} \sum_{i,\, j=1}^n \sum_{m=0}^{k-1} u_i^m u_j^{k-m-1}
\psi_{i\alpha}\mu_\beta^\alpha\psi_j^\beta \psi_{i\lambda}\mu_\nu^\lambda
\psi_j^\nu.
\eeqa
Using 
\beq
\Psi^T\Psi =\eta, ~~\Psi {\cal U}\Psi^{-1}={\rm diag}\, (u_1, \dots, u_n)
\eeq
and antisymmetry (\ref{jw15}) of the operator $\mu$ 
we obtain 
\beq
\partial_{E^k}\log\tau =
-{1\over 4}{\rm tr}\,\left(\mu \, (\mu\, {\cal U}^{k-1}
+{\cal U}\mu\, {\cal U}^{k-2} + \dots + {\cal U}^{k-1}\mu)\right).
\eeq
The next step is to derive the formula for the derivatives of the matrix
$\Psi=(\psi_{i\alpha})$. Here we use the following equations (see \cite{D3})
\beq
\partial_i \Psi = V_i \Psi
\eeq
where the matrix $V_i$ was defined in (\ref{zh26}). From this it follows that
\beq\label{bo13}
\partial_{E^k} \Psi = V^{(k)}\Psi
\eeq
with the matrix $V^{(k)}$ of the form
\beq
V^{(k)}_{ij}={u^k_i-u^k_j\over u_i-u_j} V_{ij}.
\eeq
Doing the calculations similar to the above and using the fact that
the (co)vector (\ref{bo10})
\beq\label{bo12}
H_\alpha=c^\nu_{\nu\al}=\sum_{i=1}^n {\psi_{i\alpha}\over \psi_{i1}}
\eeq
we easily obtain the second line in the formula (\ref{bo9}).
\epp

\medskip
\noindent
{\bf Remark 4.3.} For a semisimple Frobenius manifold one can uniquely
reconstruct the first derivatives $G$-function from the Frobenius
structure solving
the system (\ref{bo8}), (\ref{bo9}) for $k=2$, \dots, $n-1$ 
together with $\partial_e G=0$.
Indeed, in the canonical coordinates the coefficients of the linear system
is the Vandermonde matrix of $u_1$, \dots, $u_n$. Particularly, for any
smooth projective variety $X$ with semisimple quantum cohomology
this give a practical way to express elliptic Gromov - Witten invariants
via rational ones (see examples in Concluding Remarks below).

\medskip
We now compute derivatives of the second part
\beq
F^{(1)} ={1\over 24} \log\det M^\alpha_\beta(t,t_X)
\eeq
in the genus 1 free energy (\ref{zh19}). We will use the following formula
for this function (see \cite{DZ2})
\beq\label{bo11}
F^{(1)} ={1\over 24} \left(\log \prod_{i=1}^n \sigma_i -
\log \prod_{i=1}^n \psi_{i1} \right)
\eeq
where
\beq
\sigma_i =\psi_{i\alpha} t^\alpha_X.
\eeq
It follows that
\beq
{\partial F^{(1)}\over \partial t^\alpha_X} ={1\over 24} \sum_{i=1}^n
{\psi_{i\alpha}\over \sigma_i}.
\eeq

Let us define the operator of Lie derivative of a function
$F^{(1)}(t,t_X)$ along a vector field $v=v^\alpha(t)\partial_\alpha$ on
$M$
\beq
Lie_v F^{(1)}(t,t_X) = v^\alpha {\partial F^{(1)}\over \partial t^\alpha}
+ \left( v^\alpha\right)_X {\partial F^{(1)}\over \partial t^\alpha_X}.
\eeq
\smallskip

\begin{lemma} The  following formulae hold true for the derivatives
of the function (\ref{bo11})
\eqa
Lie_e F^{(1)}(t,t_X) &=& 0\label{bo14}\\
Lie_E  F^{(1)}(t,t_X) &=&{n\over 24}\label{bo15}\\
Lie_{E^k}  F^{(1)}(t,t_X) &=& {k\over 24} \left< E^{k-1},H\right>
\label{bo16}\\
t^\alpha_X {\partial  F^{(1)}(t,t_X)\over \partial t^\alpha_X} &=& {n\over
24}.\label{bo17}
\eeqa
Here $H$ is defined in (\ref{bo10}).
\end{lemma}

Proof easily follows from  (\ref{bo13}), (\ref{bo12}).

\medskip
We proceed now to the derivation of the genus one Virasoro constraints.
First, the $L_{-1}$ Virasoro constraint reads
\eqa
&&{\cal A}_{-1,1}\equiv 
\sum_{p\ge 1}
\tilde T^{\alpha,p} {\partial {\cal F}_1\over \partial T^{\alpha,p-1}}
\nn\\
&&
=\sum_{p\ge 1} 
\tilde T^{\alpha,p} {\partial t^\gamma\over \partial T^{\alpha,p-1}}
{\partial {\cal F}_1\over \partial t^\gamma}
+\sum_{p\ge 1} 
\tilde T^{\alpha,p} \partial_X\left( {\partial t^\gamma\over
\partial T^{\alpha, p-1}}\right) {\partial {\cal F}_1\over \partial
t_X^\gamma}
\nn\\
&&
=-\delta_1^\gamma {\partial {\cal F}_1\over \partial t^\gamma}=-Lie_e
{\cal F}_1=0.
\eeqa
In this computation we used vanishing of $S_{-1}$ on the solution
$t^{(0)}$.

Let us now prove validity of the genus 1 $L_0$ Virasoro constraint
\eqa
&&{\cal A}_{0,1} \equiv 
\left[ \sum_{p\ge 0} \left< \left( p+\mu+{1\over 2}\right)\tilde T^p,
\partial_{T^p}\right>
+ \sum_{p\geq 1}\sum_{1\le r\le p} 
\left< R_r \tilde T^p, \partial_{T^{p-r}}\right>
\right] {\cal F}_1+{1\over 4} {\rm tr}\,\left( {1\over 4}-\mu^2\right)
\nn\\
&&
= {\partial{\cal F}_1\over \partial t^\gamma}
\left[ \sum_{p\ge 0} \left< \left( p+\mu+{1\over 2}\right)\tilde T^p,
\partial_{T^p}\right>
+ \sum_{p\geq 1}\sum_{1\le r\le p}
 \left< R_r \tilde T^p, \partial_{T^{p-r}}\right>
\right]\, t^\gamma  
\nn\\
&&
+{\partial{\cal F}_1\over \partial t^\gamma_X}
\left[ \sum_{p\ge 0} \left< \left( p+\mu+{1\over 2}\right)\tilde T^p,
\partial_{T^p}\right>
+ \sum_{p\geq 1}\sum_{1\le r\le p}
 \left< R_r \tilde T^p, \partial_{T^{p-r}}\right>
\right]\, t^\gamma _X
+{1\over 4} {\rm tr}\,\left( {1\over 4}-\mu^2\right)
\nn\\
&&
=-Lie_E {\cal F}_1 -\left( {1\over 2}+\mu\right)^\nu_1 \,
\frac{\pal t^\gamma}{\pal T^{\nu,0}}\, {\partial{\cal
F}_1\over \partial t^\gamma_X}+{1\over 4} {\rm tr}\,\left( {1\over
4}-\mu^2\right)
\nn\\
&&
=-{n\over 24} - {n\,d\over 48} +{1\over 4}{\rm tr}\, \mu^2 -{n\over 24} 
\left({1-d\over 2}\right) +{n\over 16} -{1\over 4}{\rm tr}\, \mu^2=0.
\eeqa
In this computation we used vanishing of the symmetry $S_0$ on the
solution $t^{(0)}$, and also the formulae 
(\ref{bo8}), (\ref{bo15}), (\ref{bo17}) and the fact that 
$\mu^\nu_1=-\frac{d}2\,\delta^\nu_1$.

In a similar, although more involved way, using vanishing of the
symmetries $S_1$ and $S_2$ on the solution $t^{(0)}$ one can prove 
validity of the genus one $L_1$ and $L_2$ Virasoro constraints.

Now, since we proved the Virasoro constraints $L_m$, up to the genus 1,
for $-1\leq m \leq 2$, from the commutation relation
$$
[L_m, L_1]=(m-1) L_{m+1}
$$
we derive that these constraints hold true also for any $m\geq -1$.
Main Theorem is proved.
 
\setcounter{equation}{0}
\bigskip
\section{Concluding remarks}\par  
\medskip
\noindent
1. We consider the results of this paper as a strong support of the
conjectural relation between semisimple Frobenius manifolds and
integrable hierarchies of the KdV type (see Introduction above).
It gives also a practical algorithm to reconstruct the integrable
hierarchy starting from a given semisimple Frobenius manifold. The
algorithm is based on some more strong requirement that for the
tau-function $\tau(T)$ of an {\it arbitrary} solution
$t_\alpha=t_\alpha(T)$ of the hierarchy the function 
$$
\tau(T) +\delta \tau(T):= \tau(T) + \sigma L_m\tau(T)+O(\sigma^2)
$$
for any $m$ is again the tau-function of another solution of the
hierarchy in the linear approximation in the small parameter $\sigma$.
Recall that, by definition of the tau-function
$$
t_\alpha(T) = {\partial^2\log\tau(T)\over \partial T^{\alpha,0} \partial
X}, ~~\alpha=1, \dots, n.
$$
In other words, we postulate that our Virasoro operators correspond to the
symmetries 
$$
t_\alpha\mapsto t_\alpha +\sigma S_{m}[t]_\al +O(\sigma^2), 
\quad S_{m}[t]_\al:=
\frac{\partial^2}{\partial X \partial T^{\alpha,0}}\left(\frac{
L_m\tau}{\tau}\right)
$$
of the hierarchy
$$
\partial_{T^{\alpha,p}}t =\partial_X K^{(0)}_{\al,p}(t) +
\epsilon^2 \partial_X
K^{(1)}_{\al,p}(t, t_X, t_{XX}) 
+\epsilon^4 \partial_X K^{(2)}_{\al,p}(t, t_X, t_{XX},
t_{XXX}, t^{IV})+\dots
$$
in all orders in $\epsilon$. This algorithm will give us a recursion
procedure for computing of the coefficients $K^{(r)}(t, t_X, \dots,
t^{(2r)})$ for any $r$. The terms $K^{(0)}$ and $K^{(1)}$ of the hierarchy
have already been constructed in \cite{D1,D3} and \cite{DZ2} 
resp. We are going to study
the structure of higher terms in a subsequent publication.

\medskip
\noindent
2. As we have already mentioned (see Remark 4.3 above) the formulae
(\ref{bo7}), (\ref{bo8}), (\ref{bo9})
 give a simple way to compute the elliptic Gromov - Witten invariants
of those smooth projective varieties $X$ for which the quantum cohomology
is semisimple. We give here two examples of application of this
method.
The equations (\ref{bo7}), (\ref{bo8}), (\ref{bo9}) 
for the derivatives of the $G$-function can be recast
into an elegant form using a generating function
\eqa
&&\partial_{[e-zE]^{-1}} G = \sum_{k=0}^\infty z^k\partial_{E^k} G
\nn\\
&&
={z\over 24} \left( \left< \mu \left( {H\over e-z E}\right), {e\over e-z
E} \right> - 6\, {\rm tr}\, \left( \mu\, (1-z {\cal U})^{-1}\right)^2
\right)
\eeqa
Here $z$ is an indeterminate. 

\vskip 0.3cm
{\bf Example 5.1 \ }{\it $G$-function for quantum cohomology on 
$CP^1\times CP^1$.} \par

The primary free energy for $CP^1\times CP^1$ is given by \cite{DI}
\beq
F=\frac12\, (t^1)^2\, t^4+t^1\,t^2\,t^3+(t^4)^{-1}\,f(z_1,z_2),
\eeq
here 
\beq
z_1=t^2+2\,\log(t^4),\quad z_2=t^3+2\,\log(t^4),
\eeq
and
\beq
f(z_1,z_2)=\sum_{k,l\ge 0,k+l\ge 1}\frac{N^{(0)}_{k,l}}{(2(k+l)-1)!}\,
e^{k\, z_1+l\,z_2}.
\eeq
where $N^{(0)}_{k,l}$ are the numbers of rational curves on $CP^1\times CP^1$
with bidegree $(k,l)$ which pass through $2(k+l)-1$ points. They are 
defined recursively by the following formula
with the initial condition $N^{(0)}_{1,0}=1$:
$$
N^{(0)}_{k,l}=\sum (k_1\,l_2+k_2\,l_1)\,l_2\,\left(k_1\,\left(
\begin{array}{c} 2(k+l)-4\\ 2(k_1+l_1)-2\end{array}\right)
-k_2\,\left(\begin{array}{c} 2(k+l)-4\\ 2(k_1+l_1)-3\end{array}\right)\right)
N^{(0)}_{k_1,l_1}N^{(0)}_{k_2,l_2},
$$ 
here the summation is taken over $k_i, l_i$ with
 $k_1+k_2=k, l_1+l_2=l, k_i\ge 0,
l_i\ge 0$.

The derivatives $G_i := \partial G / \partial t^i$ of the $G$-function for 
$CP^1 \times CP^1$ are determined from the system (recall that $G_1$ is
always equal to 0)
\eqa
&&2 G_2 + 2 G_3 - t^4 G_4 =-{1\over 3},
\nn\\
&&
\left[ f_2- f_{22}- f_{12}\right] G_2 
+\left[ f_1- f_{12}- f_{11}\right] G_3
+4\,t^4 G_4
= {1\over 12}\left({{f_{122}}} + {{f_{112}}}\right),
\nn\\ 
&&
\left[ 4\,f_2- 12\,f_{22}- 4\,f_1- 
f_{22}\,f_1
- 8\,f_{12} 
- f_2\,f_{12}+ 2\,f_{22}\,f_{12}+ {{f_{12}}^2}
+ 4\,f_{11}+f_{22}\,f_{11}\right] G_2
\nn\\
&&
+\left[ - 4\,f_2+ 4\,f_{22}+4\,f_1- 8\,f_{12}
-f_1\,f_{12}+ {{f_{12}}^2} 
- 12\,f_{11}
- f_2\,f_{11}+ f_{22}\,f_{11}+ 2\,f_{12}\,f_{11}\right] G_3
\nn\\
&&
={1\over 3} \left({{2\,f_{22}}} 
+ {{4\,f_{12}}}+{{2\,f_{11}}}\right)
-{1\over 12}\left( {{f_{12}\,f_{122}}} 
+ {{f_{122}\,f_{11}}} + {{f_{22}\,f_{112}}} 
+ {{f_{12}\,f_{112}}}\right).\nn
\eeqa

It turns out that the G-function has the form 
\beq
G=-\frac1{12}\,t^2-\frac1{12}\,t^3+
\sum_{k,l\ge 1}\frac{ N^{(1)}_{k,l}}{(2k+2l)!}\,e^{k\,z_1+l\,z_2}
\eeq
where  $N^{(1)}_{k,l}=N^{(1)}_{l,k}$ 
are constants. We list in Table 1  the numbers $N^{(0)}_{k,l}$ and
$N^{(1)}_{k,l}$ with $k+l\le 14$. 
For $1\leq k \leq 14$ we have
\eqa
&&N^{(0)}_{k,0}=N^{(0)}_{0,k}=\delta_{k,1},\quad 
N^{(0)}_{k,1}=N^{(0)}_{1,k}=1,
\nn\\
&&N^{(1)}_{k,1}=N^{(1)}_{1,k}=0.\nn
\eeqa
This agrees with the definition of GW invariants (for any $k$).
The remaining numbers are put in the table below. 
Comparison with the results of \cite{Vakil} suggests that 
the numbers $N^{(1)}_{k,l}$ coincide with the
numbers of irreducible elliptic curves in the class $k S+l F$ on the
rational ruled surface $F_0 \simeq {\bf CP}^1 \times {\bf CP}^1$
(in the
notation of 
\cite{Vakil}, see Table 1 there). 
\par

\begin{table}

\caption{List of some numbers $N^{(0)}_{k,l}$ and $N^{(1)}_{k,l}$ 
for $CP^1\times CP^1$}
\begin{center}
\begin{tabular}{|l|l|l|}
\hline
$(k,l)$ & $N^{(0)}_{k,l}$  & $N^{(1)}_{k,l}$ \\ \hline
$(2,2)$ & $12$             & $1$                \\ \hline
$(3,2)$ & $96$             & $20$                \\ \hline
$(4,2)$ & $640 $             & $240 $                \\ \hline
$(3,3)$ & $3510 $             & $1920 $                \\ \hline
$(5,2)$ & $3840 $             & $2240 $                \\ \hline
$(4,3)$ & $87544 $             & $87612 $                \\ \hline
$(6,2)$ & $21504 $             & $17920 $                \\ \hline
$(5,3)$ & $1763415 $             & $2763840 $                \\ \hline
$(4,4)$ & $6508640 $             & $12017160 $                \\ \hline
$(7,2)$ & $114688 $             & $129024 $                \\ \hline
$(6,3)$ & $30940512 $             & $69488120 $                \\ \hline
$(5,4)$ & $348005120 $             & $1009712640 $                \\ \hline
$(8,2)$ & $589824 $             & $860160 $                \\ \hline
$(7,3)$ & $492675292 $             & $1495782720 $               
 \\ \hline
$(6,4)$ & $15090252800 $             & $62820007680 $              
  \\ \hline
$(5,5)$ & $43628131782 $             & $199215950976 $              
  \\ \hline
$(9,2)$ & $2949120 $             & $5406720 $               
 \\ \hline
$(8,3)$ & $7299248880 $             & $28742077000 $               
 \\ \hline
$(7,4)$ & $565476495360 $             & $3183404098560 $               
 \\ \hline
$(6,5)$ & $4114504926336 $             & $26965003723840 $               
 \\ \hline
$(10,2)$ & $14417920 $             & $32440320 $               
 \\ \hline
$(9,3)$ & $102276100605 $             & $506333947840 $               
 \\ \hline
$(8,4)$ & $19021741768704 $             & $138871679557632 $               
 \\ \hline
$(7,5)$ & $318794127432450 $             & $2824624505793600 $               
 \\ \hline
$(6,6)$ & $780252921765888 $             & $7337244206710400 $               
 \\ \hline

$(11,2)$ & $69206016 $             & $187432960$               
 \\ \hline
$(10,3)$ & $1370760207040 $             &  $8327258171820 $              
 \\ \hline
$(9,4)$ & $588743395737600 $             & $5402199925555200 $               
 \\ \hline
$(8,5)$ & $21377025195016320 $             & $245508475513868160 $   
              \\ \hline
$(7,6)$ & $115340307031443456 $             & $1465539494120378880 $        
        \\ \hline
$(12,2)$ & $327155712 $             & $1049624576 $               
 \\ \hline
$(11,3)$ & $17716885497906 $             &  $129517853380160 $              
 \\ \hline
$(10,4)$ & $17053897886924800 $             & $191937248700825600 $         
       \\ \hline
$(9,5)$ & $1282815980041107375 $             & $18505625758298112000 $      
          \\ \hline
$(8,6)$ & $14211230949697683456 $             & $233887641913890478080 $    
            \\ \hline
$(7,7)$ & $30814236194426422332 $             &$528646007400035492736 $     
            \\ \hline
\end{tabular}
\end{center}
\end{table}

\vskip 0.4cm
{\bf Example 5.2}\ {\it $G$-function for quantum cohomology 
on $CP^3$.}\par

The primary free energy for $CP^3$ is given by \cite{DI}
\beq
F=\frac12\,(t^1)^2\, t^4+t^1\,t^2\, t^3+\frac16\, (t^2)^3+f(z_1,z_2),
\eeq
where 
\beq
z_1=\frac{t^4}{(t^3)^2},\quad z_2=t^2+4\, \log(t^3)
\eeq
and
\beq
f(z_1,z_2)=\sum_{k\ge 1}\,\sum_{0\le l\le 2\,k}\frac{N^{(0)}_{4k-2l,l}}
{(4\,k-2\,l)!\,l!}\,z_1^l\,e^{k\,z_2}.
\eeq
The numbers $N^{(0)}_{4k-2l,l}$ are the numbers of rational 
curves of degree $k$ which
pass through $l$ points  and $4k-2l$ lines in general position
 on $CP^3$.

The derivatives of the $G$-function for $CP^3$ is determined by the system
\eqa
&&4 G_2-t_3 G_3 -2 t_4 G_4 =-1,
\nn\\
&&
\left[4\,f_{22}  -  2\, f_2-
    {{{z_1}}^2}\,f_{11}\right] G_2 +
   t_3\,\left[ -4 + f_{22}  + {1\over 2}
    {{{z_1}\,f_{12}}}\right] G_3 +
2\,{{t_3}^2}G_4
\nn\\
&&\quad =
{{f_{22}- 2\,f_{222}}\over 
      6} + {{{{{z_1}}^2}\,f_{112}}\over {12}}
\nn\\
&&
\left[ 16\,f_2- 64\,f_{22}- 
    8\,f_2\,f_{22}+ 
    24\,{{f_{22}}^2} - 
    9\,f_1- 
    24\,{z_1}\,f_1+ 
    6\,{z_1}\,f_{22}\,f_1+ 
    14\,f_{12}+64\,{z_1}\,f_{12}\right.
\nn\\
&&\left. - 
    12\,{z_1}\,f_{22}\,
     f_{12}+ 
    3\,{{{z_1}}^2}\,f_1\,
     f_{12}- 
    8\,{{{z_1}}^2}\,{{f_{12}}^2} + 
    {z_1}\,\left( -5 - 16\,{z_1} + 
       2\,{z_1}\,f_{22}+ 
       3\,{{{z_1}}^2}\,f_{12}\right) \,
     f_{11}\right] G_2
\nn\\
&&
+ {t_3\over 4}\left[
   {{ -16\,f_2+ 
          8\,{{f_{22}}^2} + 
          3\,f_1+ 6\,f_{12}+ 
          16\,{z_1}\,f_{12}- 
          2\,{{{z_1}}^2}\,{{f_{12}}^2} }} 
- {{{z_1}\,\left( 1 + 8\,{z_1} \right) \,
        f_{11}}}\right] G_3
\nn\\
&&
+{{t_3}^2}\,\left[ -8 + f_2  +{3\over 2} 
    {{{z_1}\,f_1}}\right] G_4
\nn\\
&&\quad
= {1\over 24}\left[{-12\,f_2+ 72\,f_{22}+ 
        16\,{{f_{22}}^2} - 
        96\,f_{222}- 
        12\,f_2\,f_{222}+ 
        16\,f_{22}\,f_{222}+ 
        24\,f_1+ 
        18\,{z_1}\,f_1}\right.
\nn\\
&&\quad+ 
        12\,{z_1}\,f_{222}\,
         f_1+ 16\,f_{12}- 
        84\,{z_1}\,f_{12}- 
        12\,{z_1}\,f_{22}\,
         f_{12}- 
        28\,{z_1}\,f_{222}\,
         f_{12}
\nn\\
&&\quad+ 
        96\,{z_1}\,f_{122}+ 
        6\,{z_1}\,f_2\,
         f_{122}- 
        12\,{z_1}\,f_{22}\,
         f_{122}- 
        6\,{{{z_1}}^2}\,f_1\,
         f_{122}
\nn\\
&&\quad+ 
        20\,{{{z_1}}^2}\,f_{12}\,
         f_{122}+ 
        3\,{z_1}\,f_{11}+ 
        12\,{{{z_1}}^2}\,f_{11}+ 
        6\,{{{z_1}}^2}\,f_{222}\,
         f_{11}- 
        3\,{{{z_1}}^3}\,f_{122}\,
         f_{11} 
\nn\\
&&\quad \left. - 
    {{{{z_1}}^2}\,\left( 24 - 2\,f_{22}+ 
          3\,{z_1}\,f_{12}\right) \,
        f_{112}}\right] .\nn
\eeqa

From this system it follows that
the G-function has the form
\beq
G=-\frac14\,t^2+\sum_{k\ge 1}\,\sum_{0\le l\le 2\,k}\frac{N^{(1)}_{4k-2l,l}}
{(4\,k-2\,l)!\,l!}\,z_1^l\,e^{k\,z_2},
\eeq
We can determine the numbers $N^{(1)}_{k,l}$ recursively from the
expression of the gradient $\frac{\pal G}{\pal t^2}$ in terms
of the primary free energy. 
The numbers $N^{(1)}_{4k-2l,l}+\frac{2k-1}{12}\,N^{(0)}_{4k-2l,l}$ 
represent the numbers of
elliptic curves of degree $k$ passing through $l$ points and
$4k-2l$ lines in general position as  was shown by Getzler and
Pandharipande \cite{Ge}.   

\end{document}